\documentclass[journal]{IEEEtran}
\usepackage[noadjust]{cite}

\ifCLASSINFOpdf
\else
\fi
\usepackage[cmex10]{amsmath}
\usepackage{amsthm,amssymb}
\newtheorem{theorem}{Theorem}[section]
\newtheorem{lemma}[theorem]{Lemma}
\newtheorem{corollary}[theorem]{Corollary}
\newtheorem{remark}{Remark}[section]
\newtheorem{example}[theorem]{Example}

\newcommand{\E}{{\mathbb E}}
\newcommand{\R}{{\mathbb R}}

\newcommand{\ac}{{\mathcal{A}}}

\newcommand{\K}{{\mathcal{K}}}

\newcommand{\samp}{{\mathcal{X}}}

\begin{document}
%
\title{Lower bounds for the minimax risk using $f$-divergences, and applications}
%
%
%

\author{Adityanand~Guntuboyina
\thanks{After acceptance of this manuscript, Professor Alexander
  Gushchin pointed out that Theorem~\ref{maha} appears in his
  paper~\cite{Gushchin}. Specifically, in a different notation, 
inequality~\eqref{maha.eq1} appears as Theorem 1 and
inequality~\eqref{maha.eq} appears in Section 4.3
in~\cite{Gushchin}. The proof of Theorem~\ref{maha} presented in
section~\ref{test} is
different from that in~\cite{Gushchin}. Also, except for
Theorem~\ref{maha} and the observation  
that Fano's inequality is a special case of Theorem~\ref{maha} (see
Example~\ref{kld}), there is no other overlap between this paper
and~\cite{Gushchin}.}
\thanks{Some of the material in this paper was presented at the IEEE
  International Symposium on Information Theory, Austin, TX, June 2010.}
\thanks{A. Guntuboyina is with the Department
of Statistics, Yale University, 24 Hillhouse Avenue, New Haven,
CT 06511, USA. e-mail: adityanand.guntuboyina@yale.edu}}

%
%

\markboth{TO APPEAR IN IEEE TRANSACTIONS ON INFORMATION THEORY, 2011}%
{TO APPEAR IN IEEE TRANSACTIONS ON INFORMATION THEORY, 2011}

%



\maketitle

\begin{abstract}
Lower bounds involving $f$-divergences between the underlying probability
measures are proved for the minimax risk in estimation
problems. Our proofs just use simple convexity facts. Special cases
and straightforward corollaries of our bounds include well known
inequalities for establishing minimax 
lower bounds such as Fano's inequality, Pinsker's inequality and
inequalities based on global entropy conditions. Two
applications are provided: a new minimax lower bound for the
reconstruction of 
convex bodies from noisy support function measurements and a different
proof of a recent minimax lower bound for the estimation of a
covariance matrix.
\end{abstract}


%
\IEEEpeerreviewmaketitle

\section{Introduction}\label{intro}
%
%
%
%
\IEEEPARstart{C}{onsider} an estimation problem in which
we want to estimate $\theta 
\in \Theta$ based on an observation $X$ from $\left\{ P_{\theta} ,
  \theta \in \Theta \right\}$ where each $P_{\theta}$ is a probability measure
on a sample space $\samp$. Suppose that estimators are allowed to
take values in $\ac \supseteq \Theta$ and
that the loss function is of the form $\ell(\rho)$ where $\rho$ is a metric on
$\ac$ and  $\ell: [0, \infty) \rightarrow [0, \infty)$ is a nondecreasing
function. The minimax risk for this problem is defined by
\begin{equation*}
R :=  \inf_{\hat{\theta}} \sup_{\theta \in \Theta}
\E_{\theta} \ell( \rho(\theta, \hat{\theta}(X)) ), 
\end{equation*}
where the infimum is over all measurable functions $\hat{\theta} :
\samp \rightarrow \ac$ and the expectation is taken under the
assumption that $X$ is distributed according to $P_{\theta}$. 

In this article, we are concerned with the problem of obtaining lower
bounds for the minimax risk $R$. Such bounds are useful in assessing the quality
of estimators for $\theta$. The standard approach to these bounds is to obtain a
reduction to the more tractable problem of bounding from below the
minimax risk of a multiple hypothesis testing problem. More specifically, one
considers a finite subset $F$ of the parameter space $\Theta$ and a
real number $\eta$ such that $\rho(\theta, \theta') \geq \eta$ for
$\theta, \theta ' \in F, \theta \neq \theta'$ and employs the
inequality $R \geq \ell(\eta/2)r$, 
where 
\begin{equation}\label{rdef}
r := \inf_T \sup_{\theta
  \in F} P_{\theta} \left\{T \neq \theta \right\},
\end{equation}
the infimum being over all estimators $T$ taking values in $F$. The
proof of this inequality relies on the triangle inequality satisfied
by the metric $\rho$ and can be found, for example, in~\cite[Page
1570, Proof of Theorem 1]{YangBarron} (Let us note, for the
convenience of the reader, that the notation employed by Yang and
Barron~\cite{YangBarron} differs from ours in that they use $d$ for
the metric $\rho$, $\epsilon_{n, d}$ for our $\eta$ and
$N_{\epsilon_n, d}$ for the finite set $F$. Also the proof 
in~\cite{YangBarron} involves a positive constant $A$ which 
can be taken to be 1 for our purposes. The constant $A$ arises because
Yang and Barron~\cite{YangBarron} do not require that $d$ is a metric
but rather require it to satisfy a weaker local triangle inequality
which involves the constant $A$.)

The next step is to note that $r$ is bounded from below by Bayes
risks. Let $w$ be a probability measure on $F$. The Bayes risk $\bar{r}_w$
corresponding to the prior $w$ is defined by
\begin{equation}\label{minbayes}
\bar{r}_w  :=  \inf_T \sum_{\theta \in F}
w_{\theta}P_{\theta} \left\{ T \neq \theta \right\},
\end{equation}
where $w_{\theta} := w \left\{\theta \right\}$ and the infimum is over
all estimators $T$ taking values in $F$. When $w$ is the discrete uniform
probability measure on $F$, we simply write $\bar{r}$ for $\bar{r}_w$. The
trivial inequality $r \geq \bar{r}_w$ implies that lower bounds for
$\bar{r}_w$ are automatically lower bounds for $r$. 

The starting point for the results described in this paper is
Theorem~\ref{maha}, which provides a lower 
bound for $\bar{r}_w$ involving $f$-divergences of the probability
measures $P_{\theta}, \theta \in 
F$. The $f$-divergences (\cite{AliSilvey, Csiszar66, Csiszar67,
  Csiszar67fdiv}) are a general
class of divergences between probability measures which include many
common divergences\slash distances like the Kullback Leibler
divergence, chi-squared divergence, total variation distance, Hellinger
distance etc. For a \textit{convex} function $f:[0, \infty) \rightarrow \R$
satisfying $f(1) = 0$, the $f$-divergence between two probabilities $P$ and
$Q$ is given by   
\begin{equation*}
D_f(P||Q)  :=  \int f\left( \frac{dP}{dQ} \right) dQ 
\end{equation*}
if $P$ is absolutely continuous with respect to $Q$ and $\infty$
otherwise.

Our proof of Theorem~\ref{maha} presented in section~\ref{test} is
extremely simple. It just relies on the convexity of the function 
$f$ and the standard result that $\bar{r}_w$ has the following exact
expression: 
\begin{equation}\label{betexp}
  \bar{r}_w = 1 - \int_{\samp} \max_{\theta \in F}\left\{
w_{\theta}p_{\theta}(x)\right\}  d\mu(x),
\end{equation}
where $p_{\theta}$ denotes the density of $P_{\theta}$ with respect to
a common dominating measure $\mu$ (for example, one can take $\mu :=
\sum_{\theta \in F}P_{\theta}$).

We show that Fano's inequality is a special case (see Example~\ref{kld}) of
Theorem~\ref{maha}, obtained by taking $f(x) = x \log x$. Fano's
inequality is used extensively in the
nonparametric statistics literature for obtaining minimax lower
bounds, important works being~\cite{IbragimovHasminskii77, IbragimovHasminskii80,
  Ibragimov:Hasminskii:81book, Hasminskii78, Birge83, Birge:86newZW,
  YangBarron}. In the special case when $F$ has only
two points, Theorem~\ref{maha} gives a sharp inequality relating the
total variation distance between two probability measures to
$f$-divergences (see Corollary~\ref{Neq2}). When $f(x)= x\log x$,
Corollary~\ref{Neq2} implies an inequality due to
Tops{\o}e~\cite{TopsoeCapDiv} from which Pinsker's inequality can be
derived. Thus Theorem~\ref{maha} can be viewed as a generalization of
both Fano's inequality and Pinsker's inequality. 

The bound given by Theorem~\ref{maha} involves the quantity $J_f :=
\inf_Q \sum_{\theta \in F} D_f(P_{\theta}||Q)/|F|$, where the infimum
is over all probability measures $Q$ and $|F|$ denotes the cardinality
of the finite set $F$. It is usually not possible to calculate $J_f$
exactly and in section~\ref{jsbounds}, we provide upper bounds for
$J_f$. The main result of this section, Theorem~\ref{ybgen}, provides
an upper bound for $J_f$ based on approximating the set of $|F|$
probability measures $\left\{P_{\theta}, \theta \in F\right\}$ by a
smaller set of probability measures. This result is motivated by and a
generalization to $f$-divergences of a result of Yang and
Barron~\cite{YangBarron} for the Kullback-Leibler divergence. 

In section~\ref{gent}, we use the inequalities proved in
sections~\ref{test} and~\ref{jsbounds} to
obtain minimax lower bounds involving only global metric entropy
attributes. Of all the lower bounds presented in this paper,
Theorem~\ref{myb.thm}, the main result of section~\ref{gent}, is the most
application-ready method. In order to apply this in a particular
situation, one only needs to determine suitable bounds on global
covering and packing numbers of the parameter space $\Theta$ and the
space of probability measures $\left\{ P_{\theta}, \theta \in \Theta
\right\}$ (see section~\ref{suppest} for an application). 

Although the main results of sections~\ref{test}
and~\ref{jsbounds} hold true for all $f$-divergences,
Theorem~\ref{myb.thm} is stated only for the Kullback-Leibler divergence,
chi-squared divergence and the divergences based on $f(x) = x^l -
1$ for $l > 1$. The reason behind this is that Theorem~\ref{myb.thm}
is intended 
for applications where it is usually the case that the  underlying
probability measures $P_{\theta}$ are product measures and divergences
such as the Kullback-Leibler divergence and chi-squared divergence can
be computed for product probability measures. 

The inequalities given by Theorem~\ref{myb.thm} for the chi-squared
divergence and divergences based on $f(x) = x^l - 1$ for $l > 1$ are
new while the inequality for the Kullback-Leibler divergence is due to
Yang and Barron~\cite{YangBarron}. There turn out to be qualitative
differences between these inequalities in the case of estimation
problems involving finite dimensional parameters where the inequality
based on chi-squared divergence gives minimax lower bounds having the
optimal rate while the one based on the Kullback-Leibler divergence
only results in sub-optimal lower bounds. We shall explain this
happening in section~\ref{gent} by means of elementary examples. 

We shall present two applications of our bounds. In
section~\ref{suppest}, we shall prove a new lower bound for the minimax
risk in the problem of estimation/reconstruction of a $d$-dimensional
convex body from noisy measurements of its support function in $n$
directions. In section~\ref{covmat}, we shall provide a different 
proof of a recent result by Cai, Zhang and
Zhou~\cite{CaiZhangZhou2009} on covariance matrix estimation. 
\section{Lower bounds for the testing risk $\bar{r}_w$}\label{test}
We shall prove a lower bound for $\bar{r}_w$ defined
in~\eqref{minbayes} in terms of $f$-divergences. We shall assume that
the $N := |F|$ probability measures $P_{\theta}, \theta \in F$ are all
dominated by a sigma finite measure $\mu$ with densities $p_{\theta},
\theta \in F$. In terms of these densities, $\bar{r}_w$ has the exact
expression given in~\eqref{betexp}. A trivial consequence
of~\eqref{betexp} that we shall often use in the sequel is that
$\bar{r} \leq 1 - 1/N$ (recall that $\bar{r}$ is $\bar{r}_w$ in the
case when $w$ is the uniform probability measure on $F$).  
\begin{theorem}\label{maha}
Let $w$ be a probability measure on $F$. Define $T: \samp \rightarrow
F$ by $T(x)  :=  \arg 
\max_{\theta \in F} \left\{w_{\theta} 
  p_{\theta}(x)\right\}$, where $w_{\theta} := w\left\{ \theta\right\}$.
For every convex function $f: [0, \infty) \rightarrow \R$ and every
probability measure $Q$ on $\samp$, we have 
\begin{equation}\label{maha.eq}
\sum_{\theta \in F}w_{\theta}D_f(P_{\theta}||Q)   \geq
Wf\left( \frac{1-\bar{r}_w}{W}\right) + (1-W)f\left(
  \frac{\bar{r}_w}{1-W}\right),
\end{equation}
where $W := \int_{\samp} w_{T(x)}dQ(x)$. In particular,
taking $w$ to be the uniform probability measure, we get that
\begin{equation}\label{maha.eq1}
\sum_{\theta \in F} D_f(P_{\theta}||Q)   \geq 
f\left( N(1 - \bar{r})\right) + (N - 1)f\left(
  \frac{N \bar{r}}{N - 1}\right).
\end{equation}
\end{theorem}
The proof of this theorem relies on a simple application of the
convexity of $f$ and it is presented below.
\begin{IEEEproof}
We may assume that all the weights $w_{\theta}$ are strictly
positive and that the probability measure $Q$ has a density $q$ with
respect to $\mu$. We start with a simple inequality for nonnegative
numbers $a_{\theta}, \theta \in F$ with $\tau := \arg \max_{\theta \in
F} \left\{w_{\theta} a_{\theta} \right\}$. We first write
\begin{equation*}
\sum_{\theta \in F} w_{\theta}f(a_{\theta})  =  w_{\tau}f(a_{\tau}) + (1-w_{\tau})\sum_{\theta \neq \tau}
\frac{w_{\theta}}{1-w_{\tau}}f(a_{\theta}) 
\end{equation*}
and then use the convexity of $f$ to obtain that the quantity
$\sum_{\theta}w_{\theta} f(a_{\theta})$ is bounded from below by
\begin{equation*}
  w_{\tau}f(a_{\tau}) + (1-w_{\tau}) f \left(
  \frac{\sum_{\theta \in F} w_{\theta} a_{\theta} -  w_{\tau}a_{\tau}}{1-w_{\tau}}\right). 
\end{equation*}
We now fix $x \in \samp$ such that $q(x) > 0$ and apply the inequality
just derived to $a_{\theta} := p_{\theta}(x)/q(x)$. Note that in this
case $\tau = T(x)$. We get that
\begin{equation}\label{aux1}
\sum_{\theta \in F} w_{\theta}f \left( \frac{p_{\theta}(x)}{q(x)} \right)
 \geq A(x) + B(x),
\end{equation}
where 
\begin{equation*}
A(x)  := w_{T(x)} f\left( \frac{p_{T(x)}(x)}{q(x)} \right) 
\end{equation*}
and
\begin{equation*}
B(x) := (1-w_{T(x)}) f \left( \frac{\sum_{\theta \in
      F}w_{\theta}p_{\theta}(x) 
  - w_{T(x)}p_{T(x)}(x)}{(1-w_{T(x)})q(x)} \right).
\end{equation*}
Integrating inequality~\eqref{aux1} with respect to the probability 
measure $Q$, we get that the term $\sum_{\theta \in F}
w_{\theta}D_f(P_{\theta}||Q)$ is bounded from below by 
\begin{equation*}
\int_{\samp}A(x)q(x)d\mu(x) + \int_{\samp}B(x)q(x)d\mu(x).
\end{equation*}
Let $Q'$ be the probability measure on $\samp$ having the density
$q'(x) := w_{T(x)}q(x)/W$ with respect to $\mu$. Clearly
\begin{equation*}
\int_{\samp}A(x)q(x)d\mu(x) = W \int_{\samp} f \left(\frac{p_{T(x)}(x)}{q(x)} \right) q'(x)
d\mu(x),
\end{equation*}
which, by Jensen's inequality, is larger than or equal to
$Wf((1-\bar{r}_w)/W)$. It follows similarly that
\begin{equation*}
\int_{\samp}B(x) q(x) d\mu(x)  \geq  (1-W) f\left(
  \frac{\bar{r}_w}{1-W} \right).
\end{equation*}
This completes the proof of inequality~\eqref{maha.eq}. When $w$ is
the uniform probability measure on the finite set $F$, it is obvious
that $W$ equals $1/N$ and this leads to inequality~\eqref{maha.eq1}. 
\end{IEEEproof}
Let us denote the function of $\bar{r}$ on the right hand side
of~\eqref{maha.eq1} by $g$:
\begin{equation}\label{gintang}
  g(a) := f\left(N(1-a)\right) + (N-1) f \left(\frac{Na}{N-1} \right).
\end{equation}
Inequality~\eqref{maha.eq1} provides an implicit lower bound for
$\bar{r}$. This is because $\bar{r} \in [0, 1-1/N]$ and $g$ is
non-increasing on $[0, 1-1/N]$ (as can be seen in the proof of the
next corollary in the case when $f$ is differentiable; if $f$ is not
differentiable, one needs to work with right and left derivatives
which exist for convex functions). 

The convexity of $f$ also implies trivially that $g$ is
convex, which can be used to convert the implicit
bound~\eqref{maha.eq1} into an explicit lower bound. This is the
content of the following corollary. We assume differentiability for
convenience; to avoid working with one-sided derivatives. 
\begin{corollary}\label{explow}
Suppose that $f: [0, \infty)$ is a differentiable convex function and
that $g$ is defined as in~\eqref{gintang}. Then, for every $a \in
[0,1-1/N]$, we have 
\begin{equation}\label{explow.eq}
r  \geq  \bar{r}  \geq  a + \frac{\inf_{Q}\sum_{\theta \in F}
  D_f(P_{\theta}||Q)  -  g(a)}{g'(a)},
\end{equation}
where the infimum is over all probability measures $Q$. 
\end{corollary}
\begin{IEEEproof}
Fix a probability measure $Q$. Inequality~\eqref{maha.eq1} says that
$\sum_{\theta \in F} D_f(P_{\theta} || Q) \geq g(\bar{r})$. The
convexity of $f$ implies that $g$ is also convex and hence, for every
$a \in [0, 1-1/N]$, we can write   
\begin{equation}\label{explow.temp}
\sum_{\theta \in F} D_f(P_{\theta}||Q) \geq g(\bar{r})  \geq  g(a) +
g'(a)(\bar{r} - a).
\end{equation}
Also, 
\begin{equation*}
\frac{g'(a)}{N}  =  f'\left(\frac{Na}{N-1}\right) - f'\left(N(1-a)\right). 
\end{equation*}
Because $g$ is convex, we have $g'(a) \leq g'(1-1/N) = 0$ for $a \leq 
1-1/N$ (this proves that $g$ is non-increasing on $[0,
1-1/N]$). Therefore, by rearranging~\eqref{explow.temp}, we 
obtain~\eqref{explow.eq}.
\end{IEEEproof} 
Let us now provide an intuitive understanding of
inequality~\eqref{maha.eq1}. When the probability measures
$P_{\theta}, \theta \in F$ are tightly packed i.e., when they are 
close to one another, it is hard to distinguish between them
(based on the observation $X$) and hence, the testing Bayes risk $\bar{r}$
will be large. On the other hand, when the probability measures are
well spread out, it is easy to distiguish between them and
therefore, $\bar{r}$ will be small. Indeed, $\bar{r}$ takes on its
maximum value of $1-1/N$ when the probability measures $P_{\theta},
\theta \in F$ are all equal to one another and it takes on its
smallest value of 0 when $\max p_{\theta} = \sum p_{\theta}$ i.e.,
when $P_{\theta}, \theta \in F$ are all mutually singular. 

Now, one way of measuring how packed/spread out the probability
measures $P_{\theta}, \theta \in F$ are is to consider the quantity
$\inf_Q \sum_{\theta \in F} D_f(P_{\theta}||Q)$, which is small when
the probabilities are tightly packed and large when they are spread
out. It is therefore reasonable to expect a connection between this
quantity and $\bar{r}$. Inequality~\eqref{maha.eq1} makes 
this connection explicit and precise. The fact that the function $g$ 
in~\eqref{gintang} is non-increasing means that when $\inf_Q
\sum_{\theta \in F} D_f(P_{\theta}||Q)$ is small, the lower bound on
$\bar{r}$ implied by~\eqref{maha.eq1} is large and when $\inf_Q
\sum_{\theta \in F} D_f(P_{\theta}||Q)$ is large, the lower bound on
$\bar{r}$ is small. 

Theorem~\ref{maha} implies the following corollary which provides
sharp inequalities between total variation distance and
$f$-divergences. The total variation distance between two probability
measures is defined as \emph{half} the $L^1$ distance between their
densities. 
\begin{corollary}\label{Neq2}
Let $P_1$ and $P_2$ be two probability measures on a space
$\samp$ with total variation distance $V$. For every convex function $f:
[0, \infty) \rightarrow \R$, we have 
\begin{equation}\label{tvf.eq}
\inf_Q \left( D_f(P_1||Q) + D_f(P_2||Q) \right)  \geq 
f\left( 1 + V \right) + f\left( 1 - V \right),
\end{equation}
where the infimum is over all probability measures $Q$. Moreover this
inequality is sharp in the sense that for every $V \in [0, 1]$, the
infimum of the left hand side of~\eqref{tvf.eq} over all probability
measures $P_1$ and $P_2$ with total variation distance $V$ equals the
right hand side of~\eqref{tvf.eq}. 
\end{corollary}
\begin{IEEEproof}
In the setting of Theorem~\ref{maha}, suppose that $F= \left\{1,
  2\right\}$ and that the two probability measures are $P_1$ and
$P_2$ with densities $p_1$ and $p_2$ respectively. Since $2\max (p_1,
p_2)$ equals $p_1+p_2+|p_1-p_2|$, it follows that $2\bar{r}$ equal
$1-V$. Inequality~\eqref{tvf.eq} is then a direct consequence 
of inequality~\eqref{maha.eq1}.  

The following example shows that~\eqref{tvf.eq} is sharp. Fix $V \in
[0, 1]$. Consider the space $\samp = \left\{1, 2 \right\}$ and define
the probabilities $P_1$ and $P_2$ by $P_1\left\{ 1 \right\} =
P_2\left\{2 \right\} = (1+V)/2$ and of course $P_1\left\{ 2
\right\} = P_2\left\{1 \right\} = (1-V)/2$. Then the total 
variation distance between $P_1$ and $P_2$ equals $V$. Also if we take
$Q$ to be the uniform probability measure $Q \left\{1 \right\} =
Q\left\{2 \right\} = 1/2$, then one sees that $D_f(P_1||Q) +
D_f(P_2||Q)$ equals $f(1+V) + f(1-V)$ which is same as the right hand
side in~\eqref{tvf.eq}.   
\end{IEEEproof}
What we have actually shown in the above proof is that
inequality~\eqref{tvf.eq} is sharp for the space $\samp = \left\{1, 2
\right\}$. However, the result holds in more general spaces as
well. For example, if the space is such that there exist two disjoint 
nonempty subsets $A_1$ and $A_2$ and two probability measures $\nu_1$
and $\nu_2$ concentrated on $A_1$ and $A_2$ respectively, then we can
define $P_1 := \nu_1(1+V)/2 + \nu_2 (1-V)/2$ and  $P_2 := \nu_1(1-V)/2
+ \nu_2 (1+V)/2$ so that $V(P_1, P_2) = V$ and~\eqref{tvf.eq} becomes
an equality (with $Q = \nu_1/2 + \nu_2/2$).

There exist many inequalities in the literature relating the
$f$-divergence of two probability measures to their total variation
distance. We refer the reader to~\cite{genpin} for the sharpest
results in this direction and for earlier
references. Inequality~\eqref{tvf.eq}, which is new, can be trivially
converted into an inequality between $D_f(P_1||P_2)$ and $V$ by taking
$Q = P_2$. The resulting inequality will not be sharp however and
hence will be inferior to the inequalities in~\cite{genpin}. As
stated, inequality~\eqref{tvf.eq} is a sharp inequality relating not
$D_f(P_1||P_2)$ but a symmetrized form of $f$-divergence between $P_1$
and $P_2$ to their total variation distance.

In the remainder of this section, we shall apply Theorem~\ref{maha} and
Corollary~\ref{Neq2} to specific $f$-divergences.
\begin{example}[Kullback-Leibler Divergence]\label{kld}\normalfont
Let $f(x) := x \log x$. Then $D_f(P||Q)$ becomes the Kullback-Leibler
divergence $D(P||Q)$ between $P$ and $Q$. The quantity $\sum_{\theta \in
  F}D(P_{\theta}||Q)$ is minimized when  $Q =
\bar{P} := (\sum_{\theta \in F}P_{\theta})/N$. This is a consequence
of the following identity which is sometimes referred to as the
\textit{compensation identity}, see for example~\cite[Page 1603]{TopsoeCapDiv}:
\begin{equation*}
\sum_{\theta \in F}D(P_{\theta}|| Q)  = \sum_{\theta \in F}
D(P_{\theta}|| \bar{P}) + N D(\bar{P} || Q). 
\end{equation*}
Using inequality~\eqref{maha.eq1} with  $Q = \bar{P} = (\sum_{\theta \in 
  F}P_{\theta})/N$, we obtain
\begin{equation*}
\frac{1}{N} \sum_{\theta \in F} D(P_{\theta}||\bar{P})  \geq
(1-\bar{r})\log(N(1-\bar{r})) + \bar{r} \log \left( \frac{N
    \bar{r}}{N-1} \right) . 
\end{equation*}
The quantity on the left hand side is known as the 
Jensen-Shannon divergence. It is also Shannon's mutual
information~\cite[Page 19]{CoverThomas} between the random parameter
$\theta$ distributed according to the uniform distribution on $F$ and
the observation $X$ whose conditional distribution given $\theta$
equals $P_{\theta}$. The above inequality is stronger
than the version of Fano's inequality commonly used in nonparametric 
statistics. It is
implicit in~\cite[Proof of Theorem 1]{HanVerdu} and is explicitly stated in
a slightly different form in~\cite[Theorem
3]{BirgeFano}. The proof in~\cite{HanVerdu} is based on the Fano's inequality from
information theory~\cite[Theorem 2.10.1]{CoverThomas}. To obtain the
usual form of Fano's inequality as used in statistics, we turn to
inequality~\eqref{explow.eq}. For $a_0 := (N-1)/(2N-1) \leq 1-1/N$ and the
function $g$ in~\eqref{gintang}, it can be checked that
\begin{equation*}
  g(a_0) = \frac{N^2}{2N-1} \log N + N \log \left( \frac{N}{2N-1}
  \right)
\end{equation*}
and $g'(a_0) = -N \log N$. Using inequality~\eqref{explow.eq} with
$a = a_0$, we get that  
\begin{equation*}
  \bar{r}  \geq  1 - \frac{\log ((2N-1)/N) + \frac{1}{N}\sum_{\theta \in
      F}D(P_{\theta}||\bar{P}) }{\log N}.
\end{equation*}
Since $\log ((2N-1)/N) \leq \log 2$, we have obtained
\begin{equation}\label{normfano}
r \geq \bar{r}  \geq  1 - \frac{\log 2 + \frac{1}{N}\sum_{\theta \in
      F}D(P_{\theta}||\bar{P}) }{\log N},
\end{equation}
which is the commonly used version of Fano's inequality.

By taking $f(x) = x\log x$ in Corollary~\ref{Neq2}, we get that 
\begin{equation*}
D(P_1||\bar{P}) + D(P_2||\bar{P}) \geq (1+V)\log (1 + V) + (1-V)\log (1 - V).
\end{equation*}
This inequality relating the Jensen-Shannon divergence between two
probability measures (also known as
capacitory discrimination) to their total variation distance is due to
Tops{\o}e~\cite[Equation (24)]{TopsoeCapDiv}. Our proof is slightly simpler 
than Tops{\o}e's. Tops{\o}e~\cite{TopsoeCapDiv} also explains how to
use this inequality to deduce Pinsker's inequality with sharp
constant: $D(P_1|| P_2) \geq 2 V^2$. Thus, Theorem~\ref{maha} can be 
considered as a generalization of both Fano's inequality and Pinsker's
inequality to $f$-divergences.
\end{example}
\begin{example}[Chi-Squared Divergence]\label{csd}\normalfont
Let $f(x) = x^2-1$. Then $D_f(P||Q)$ becomes the chi-squared
divergence $\chi^2(P|| Q) := \int p^2/q -1$. The function $g$ can be
easily seen to satisfy 
\begin{equation*}
  g(a) = \frac{N^3}{N-1} \left(1-\frac{1}{N} - a \right)^2 \geq N^2
  \left(1-\frac{1}{N} - a \right)^2. 
\end{equation*}
Because $\bar{r} \leq 1-1/N$, we can invert the inequality
$\inf_Q\sum_{\theta \in F} \chi^2(P_{\theta}||Q) \geq g(\bar{r})$ to
obtain 
\begin{equation}\label{chifano}
r  \geq  \bar{r}  \geq  1 - \frac{1}{N} - \frac{1}{\sqrt{N}}
\sqrt{\frac{\inf_Q \sum_{\theta \in F}  \chi^2(P_{\theta}|| Q)}{N}}.
\end{equation}
Also it follows from Corollary~\ref{Neq2} that for every two probability
measures $P_1$ and $P_2$, 
\begin{equation}\label{tridis}
\inf_Q \left( \chi^2(P_1||Q) + \chi^2(P_2||Q) \right)  \geq  2 V^2.
\end{equation}
The weaker inequality $\chi^2(P_1||\bar{P}) + \chi^2(P_2||\bar{P})
\geq 2V^2$ can be found in~\cite[Equation (11)]{TopsoeCapDiv}.  
\end{example}
\begin{example}[Hellinger Distance]\label{helld}\normalfont
Let $f(x) = 1-\sqrt{x}$. Then $D_f(P||Q) = 1 - \int \sqrt{pq} d\mu =
H^2(P, Q)/2$, where $H^2(P,Q) = \int (\sqrt{p} - \sqrt{q})^2 d
\mu$ is the square of the Hellinger distance between $P$ and $Q$. It
can be shown, using the Cauchy-Schwarz inequality, that 
$\sum_{\theta \in F}D_f(P_{\theta}||Q)$ is minimized when $Q$ has a
density with respect to $\mu$ that is proportional to $(\sum_{\theta
  \in F}\sqrt{p_{\theta}})^2$. Indeed if $u := \sum_{\theta \in F}
\sqrt{p_{\theta}}$, then 
\begin{align*}
  \sum_{\theta \in F} D_f(P_{\theta}||Q) = N - \int \sqrt{qu^2}d\mu
  \geq N - \sqrt{\int u^2 d\mu},
\end{align*}
by the Cauchy-Schwarz inequality with equality when $q$ is
proportional to $u^2$. The inequality~\eqref{maha.eq1} can then be
simplified to
\begin{equation}\label{helldtemp}
  \sqrt{1-\bar{r}} + \sqrt{(N-1)\bar{r}} \geq \sqrt{\frac{\int u^2 d\mu}{N}}.
\end{equation}
We now observe that
\begin{equation*}
  \int u^2 d\mu = N + \sum_{\theta \neq \theta '} \int
  \sqrt{p_{\theta} p_{\theta '}} d\mu = N^2 - \frac{1}{2} \sum_{\theta
  \neq \theta '} H^2(P_{\theta}, P_{\theta '}).
\end{equation*}
We let $h^2 := \sum_{\theta, \theta'}H^2(P_{\theta},P_{\theta'})/N^2$ so
that $\int u^2 d\mu = N^2(1-h^2/2)$. As a consequence, we have $\int
u^2 d\mu \leq N^2$. Also note that $\int u^2 d\mu \geq \int (\sum_{\theta}p_{\theta})
d\mu = N$. Therefore, the right hand side of the
inequality~\eqref{helldtemp} lies between 1 and $\sqrt{N}$. On the
other hand, it can be checked that, as a function of $\bar{r}$, the
left hand side of~\eqref{helldtemp} is strictly increasing from $1$ at
$\bar{r}=0$ to $\sqrt{N}$ at $\bar{r} = 1-1/N$. It therefore follows that
inequality~\eqref{helldtemp} is equivalent to $\bar{r} \geq \breve{r}$ where
$\breve{r} \in [0, 1-1/N]$ is the solution to the equation obtained by
replacing the inequality in~\eqref{helldtemp} with an equality. 

This equation can be solved in the usual way by squaring etc., until we
get a quadratic equation in $\bar{r}$ which can be solved resulting in
two solutions. One of the two solutions can be discarded by continuity
considerations (the solution has to be continuous in $\int u^2
d\mu/N$) and the fact that $\bar{r} \leq 1-1/N$. The other solution
equals $\breve{r}$ and is given by 
\begin{equation*}
  \breve{r} = 1-\frac{1}{N} - \frac{N-2}{N}\frac{h^2}{2} -
  \frac{\sqrt{N-1}}{N} \sqrt{h^2(2-h^2)}.
\end{equation*}
We have thus shown that 
\begin{equation*}
  r \geq \bar{r} \geq 1-\frac{1}{N} - \frac{N-2}{N}\frac{h^2}{2} -
  \frac{\sqrt{N-1}}{N} \sqrt{h^2(2-h^2)}.
\end{equation*}
In the case when $N=2$ and $F = \left\{1, 2 \right\}$, it is clear
that $h^2 =(H^2(P_1, P_2) + H^2(P_2, P_1))/4= H^2(P_1, P_2)/2$. Also
since $2 \bar{r}$ equals $1-V$, 
where $V$ denotes the total variation distance between $P_1$ and
$P_2$, the above inequality implies that for every pair of probability
measures $P_1$ and $P_2$, we have
\begin{equation*}
V  \leq  H(P_1, P_2) \sqrt{1- \frac{H^2(P_1, P_2)}{4}}. 
\end{equation*}
This inequality is usually attributed to Le Cam~\cite{LeCam:86book}. 
\end{example}
\begin{example}[Total Variation Distance]\normalfont
Let $f(x) = |x-1|/2$. Then $D_f(P||Q)$ becomes the total variation
distance between $P$ and $Q$. The function $g$ satisfies
\begin{equation*}
  g(\bar{r}) = \frac{1}{2} |N(1-\bar{r}) - 1| + \frac{N-1}{2}
  \left|\frac{N\bar{r}}{N-1} -1 \right|.
\end{equation*}
Since $\bar{r} \leq 1-1/N$, we have $N(1-\bar{r}) \geq 1$ and
$N\bar{r}/(N-1) \leq 1$ so that the above expression for $g(\bar{r})$
simplifies to $N-1-N\bar{r}$. Inequality~\eqref{maha.eq1}, therefore,
results in
\begin{equation*}
r \geq \bar{r}  \geq  1 - \frac{1}{N} - \frac{\inf_Q \sum_{\theta \in
    F}V_{\theta}}{N}.
\end{equation*}
where $V_{\theta}$ denotes the total variation distance between
$P_{\theta}$ and $Q$.
\end{example}
\begin{example}\label{lchi}\normalfont
Let $f(x) = x^l - 1$ where $l>1$. The case $l = 2$ has already
been considered in Example~\ref{csd}. The function $g$ has the
expression
\begin{equation*}
  g(\bar{r}) =  N^l(1-\bar{r})^l - N + (N-1)\left(\frac{N\bar{r}}{N-1}\right)^l.
\end{equation*}
It therefore follows that $\inf_Q \sum_{\theta \in F}
D_f(P_{\theta}||Q) \geq g(\bar{r}) \geq N^l(1-\bar{r})^l - N$ which
results in the inequality 
\begin{equation}\label{lchi.eq}
r  \geq  \bar{r}  \geq  1 - \left(\frac{1}{N^{l-1}}+ \frac{\inf_Q \sum_{\theta \in F}
    D_f(P_{\theta}||Q)}{N^l} \right)^{1/l} .
\end{equation}
 When $l = 2$, inequality~\eqref{lchi.eq} results in a bound that is weaker than
inequality~\eqref{chifano} although for large $N$, the two bounds are
almost the same. 
\end{example}
\begin{example}[\textit{Reverse} Kullback-Leibler divergence]\normalfont
Let $f(x) = -\log x$ so that $D_f(P||Q) = D(Q|| P)$. Then from
Corollary~\ref{Neq2}, we get that for every two probability measures
$P_1$ and $P_2$, 
\begin{equation*}
\inf_Q \left\{ D(Q|| P_1) + D(Q|| P_2) \right\} \geq  \log \left(
  \frac{1}{1 - V^2}\right). 
\end{equation*}
This can be rewritten to get
\begin{equation}\label{myrevpins}
V  \leq  \sqrt{1 - \exp \left(- \inf_Q \left\{ D(Q|| P_1) + D(Q|| P_2)
    \right\}\right)}. 
\end{equation}
Unlike Example~\ref{kld}, it is not true that $D(Q|| P_1) + D(Q||
P_2)$ is minimized when $Q = \bar{P}$. This is easy to see because
$D(\bar{P}, P_1) + D(\bar{P}, P_2)$ is finite only when $P_1 << P_2$
and $P_2 << P_1$. By taking $Q = P_1$ and $Q = P_2$, we get
that 
\begin{equation*}
V  \leq  \sqrt{1 - \exp \left( - \min \left( D(P_1|| P_2), D(P_2||
      P_1)\right) \right)}. 
\end{equation*}
The above inequality, which is clearly weaker than
inequality~\eqref{myrevpins}, can also be found in~\cite[Proof of
Lemma 2.6]{Tsybakovbook}.  
\end{example}
\section{Bounds for  $J_f$}\label{jsbounds}
In order to apply the minimax lower bounds of the previous section in
practical situations, we
must be able to bound the quantity $J_f := \inf_Q \sum_{\theta
\in F} D_f(P_{\theta}||Q)/N$ from above. We shall provide such bounds in
this section. It should be noted that for some functions $f$, 
it may be possible to calculate $J_f$ directly. For example, the
quantity $\inf_Q \sum_{\theta \in F}H^2(P_{\theta},Q)$ can be written
in terms of pairwise Hellinger distances (Example~\ref{helld}) and may
be calculated exactly for certain probability measures
$P_{\theta}$. This is not the case for most functions $f$ however.  

The following is a simple upper bound for $J_f$ which, in the case
when $f(x) = x \log x$ or Kullback-Leibler divergence, has been
frequently used in the literature (see for example~\cite{Birge83}
and~\cite{nemirovski2000}).  
\begin{align*}
  J_f &\leq \frac{1}{N} \sum_{\theta \in F}D_f(P_{\theta}||\bar{P})
  \\
&\leq \frac{1}{N^2} \sum_{\theta, \theta ' \in
    F}D_f(P_{\theta}||P_{\theta '}) \leq \max_{\theta, \theta ' \in F}
  D_f(P_{\theta}||P_{\theta '}).
\end{align*}
We observed in section~\ref{test} that $J_f$ measures the
\textit{spread} of the probability measures $P_{\theta}, \theta \in F$
i.e., how tightly packed/spread out they are. It should be clear that
the simple bound $\max_{\theta, \theta '}D_f(P_{\theta}||P_{\theta
  '})$ does not adequately describe this aspect of $P_{\theta}, \theta
\in F$ and it is therefore desirable to look for alternative upper
bounds for $J_f$ that capture the notion of spread in a better way. 

In the case of the Kullback-Leibler divergence, Yang and
Barron~\cite[Page 1571]{YangBarron} provided such an upper bound for
$J_f$. They showed that for any finite set $\left\{Q_{\alpha}: \alpha
  \in G \right\}$ of probability measures,
\begin{equation}\label{jsyb}
  \frac{1}{N}  \sum_{\theta \in F} D(P_{\theta}||\bar{P})  \leq
 \log |G| + \max_{\theta \in F} \min_{\alpha \in
  G}D(P_{\theta}||Q_{\alpha}).
\end{equation}
Let us now take a closer look at this beautiful inequality of Yang and
Barron~\cite{YangBarron}. The $|G|$ probability measures $Q_{\alpha},
\alpha \in G$ can be viewed as an approximation of the $N$ probability
measures $P_{\theta}, \theta \in F$. The term $\max_{\theta}
\min_{\alpha}D(P_{\theta}||Q_{\alpha})$ then denotes the approximation
error, measured via the Kullback-Leibler divergence. The right hand
side of inequality~\eqref{jsyb} can therefore be made small if it is
possible to choose not too many probability measures $Q_{\alpha}$ which
well approximate the given set of probability measures $P_{\theta}$. 

It should be clear how the upper bound~\eqref{jsyb} measures the spread of
the probability measures $P_{\theta}, \theta \in F$. If the
probabilities are tightly packed, it is possible to approximate them
well with a smaller set of probabilities and then the bound will be
small. On the other hand, if $P_{\theta}, \theta \in F$ are well
spread out, we need more probability measures for approximation and
consequently the bound will be large. 

Another important aspect of inequality~\eqref{jsyb} is that it can be
used to obtain lower bounds for $R$ depending only on global metric
entropy properties of the parameter space $\Theta$ and the space of
probability measures $\left\{P_{\theta}, \theta \in \Theta
\right\}$ (see section~\ref{gent}). On the other
hand, the evaluation of inequalities resulting  
from the use of $J_f \leq \max_{\theta,
  \theta '} D(P_{\theta}||P_{\theta '})$ requires knowledge of both metric
entropy and the existence of certain special localized subsets. We
refer the reader to~\cite{YangBarron} for a detailed
discussion of these issues. 

The goal of this section is to generalize inequality~\eqref{jsyb}
to $f$-divergences. The main result is given below. In
section~\ref{gent}, we shall use this theorem along with the results
of the previous section to come up with minimax lower bounds involving
global entropy properties.  
\begin{theorem}\label{ybgen}
Let $Q_{\alpha}, \alpha \in G$ be $M := |G|$ probability measures
having densities $q_{\alpha}, \alpha \in G$ with respect to $\mu$ and
let $j:F \rightarrow G$ be a mapping from $F$ to $G$. For every convex
function $f: [0, \infty) \rightarrow \R$, we have 
\begin{equation}\label{ybgen.eq}
J_f \leq 
\frac{1}{N}\sum_{\theta \in F} \int_{\samp} \frac{q_{j(\theta)}}{M} f \left(
  \frac{Mp_{\theta}}{q_{j(\theta)}} \right) d\mu +  \left( 1 - \frac{1}{M}
  \right) f(0).
\end{equation}
\end{theorem}
\begin{IEEEproof}
Let $\bar{Q} := \sum_{\alpha \in G} Q_{\alpha}/M$ and $\bar{q} :=
\sum_{\alpha \in G} q_{\alpha}/M$. Clearly for each $\theta \in F$, we
have
\begin{equation*}
D_f(P_{\theta}||\bar{Q}) = \int_{\samp} \bar{q} \left[f \left( \frac{p_{\theta}}{\bar{q}}\right)
- f(0)\right] d\mu + f(0).
\end{equation*}
The convexity of $f$ implies that the map $y
\mapsto y[f(a/y) - f(0)]$ is non-increasing for every nonnegative
$a$. Using this and the fact that $\bar{q} \geq q_{j(\theta)}/M$, we
get that for every $\theta \in F$, 
\begin{equation*}
D_f(P_{\theta}||\bar{Q})  \leq 
\int_{\samp} \frac{q_{j(\theta)}}{M} \left[ f \left(
  \frac{Mp_{\theta}}{q_{j(\theta)}} \right) - f(0) \right]d\mu + f(0).
\end{equation*}
Inequality~\eqref{ybgen.eq} now follows as a consequence of the inequality $J_f \leq \sum_{\theta \in
  F} D_f(P_{\theta}||\bar{Q})/N$.
\end{IEEEproof}
In the following examples, we shall demonstrate that Theorem~\ref{ybgen}
is indeed a generalization of the bound~\eqref{jsyb} to
$f$-divergences. We shall also see that Theorem~\ref{ybgen}
results in inequalities that have the same qualitative structure
as~\eqref{jsyb}, at least for the convex functions $f$ of interest such as
$x^l-1, l > 1$ and $(\sqrt{x} - 1)^2$.
\begin{example}[Kullback-Leibler divergence]\normalfont
Let $f(x) = x\log x$. In this case, $J_f$ equals $\sum_{\theta \in F}
D(P_{\theta}||\bar{P})/N$ and invoking inequality~\eqref{ybgen.eq}, we get that
\begin{equation*}
\frac{1}{N}  \sum_{\theta \in F} D(P_{\theta}||\bar{P})  \leq  \log M +
\frac{1}{N} \sum_{\theta \in F} D(P_{\theta}||Q_{j(\theta)}).
\end{equation*}
Inequality~\eqref{jsyb} now follows if we choose $j(\theta) := \arg \min_{\alpha \in
G} D(P_{\theta}||Q_{\alpha})$. Hence Theorem~\ref{ybgen} is indeed a
generalization of~\eqref{jsyb}.
\end{example}
\begin{example}\label{chijsd}\normalfont
Let $f(x) = x^l -1$ for $l > 1$. Applying
inequality~\eqref{ybgen.eq}, we get that 
\begin{equation*}
J_f \leq M^{l-1}\left( \frac{1}{N} \sum_{\theta \in F}
D_f(P_{\theta}||Q_{j(\theta)}) + 1 \right) - 1.
\end{equation*}
By choosing $j(\theta) = \arg \min_{\alpha \in G}
D_f(P_{\theta}||Q_{\alpha})$, we get that
\begin{equation}\label{ljsd.eq}
J_f  \leq M^{l-1} \left( \max_{\theta \in F} \min_{\alpha \in G}
D_f(P_{\theta}||Q_{\alpha}) + 1 \right)  - 1.
\end{equation}
In particular, in the case of the chi-squared divergence i.e., when
$l=2$, the quantity $J_f = \inf_Q \sum_{\theta \in F}
\chi^2(P_{\theta}||Q)/N$ is bounded from above by
\begin{equation}\label{chijsd.eq}
M \left( \max_{\theta \in F} \min_{\alpha \in
    G}\chi^2(P_{\theta}||Q_{\alpha}) + 1 \right) - 1.
\end{equation}
Just like~\eqref{jsyb}, each of the above two inequalities is also a
function of the number of probability measures $Q_{\alpha}$ and the
approximation error which is now measured in terms of the chi-squared
divergence.  
\end{example}
\begin{example}[Hellinger distance]\normalfont
Let $f(x) = (\sqrt{x} - 1)^2$ so that $D_f(P|| Q) = H^2(P, Q)$, the
square of the Hellinger distance between $P$ and $Q$. Using
inequality~\eqref{ybgen.eq}, we get that
\begin{equation*}
J_f\leq 2 - \frac{1}{\sqrt{M}} \left(2 - \frac{1}{N} \sum_{\theta \in
    F} H^2(P_{\theta}, Q_{j(\theta)}) \right).
\end{equation*}
If we now choose $j(\theta) := \arg \min_{\alpha \in G}
H^2(P_{\theta}, Q_{\alpha})$, then we get
\begin{equation*}
J_f \leq 2 - \frac{1}{\sqrt{M}} \left(2 - \max_{\theta \in F}
    \min_{\alpha \in G} H^2(P_{\theta}, Q_{\alpha}) \right).
\end{equation*}
Notice, once again, the trade-off between $M$ and the
approximation error which is measured in terms of the Hellinger
distance.
\end{example}
\section{Bounds involving global entropy}\label{gent}
In this section, we shall apply the results of the previous two
sections to obtain lower bounds for the minimax risk $R$ depending
only on global metric entropy properties of the parameter space. The
theorem is stated below, but we shall need to establish some notation
first.  
\begin{enumerate}
\item  For $\eta > 0$, let $N(\eta) \geq 1$ be a real number for
which there exists a finite subset $F \subseteq \Theta$ with
cardinality $\geq N(\eta)$ satisfying  $\rho(\theta, \theta ') \geq \eta$ whenever
$\theta, \theta ' \in F$ and $\theta \neq \theta '$. In other words,
$N(\eta)$ is a lower bound on the $\eta$-packing number of the metric
space $(\Theta, \rho)$.
\item For a convex function $f: [0, \infty) \rightarrow \R$ satisfying $f(1)=0$, a subset $S
  \subseteq \Theta$ and a positive real number $\epsilon$, let 
$M_f(\epsilon, S)$ be a positive real number for which there
exists a finite set $G$ with cardinality $\leq M_f(\epsilon, S)$
and probability measures $Q_{\alpha}, \alpha \in G$ such that
$\sup_{\theta \in S} \min_{\alpha \in G} D_f(P_{\theta}||Q_{\alpha})
\leq \epsilon^2$. In other words, $M_{f}(\epsilon, S)$ is an upper bound on the
$\epsilon$-covering number of the space $\left\{P_{\theta}: \theta \in
S \right\}$ when distances are measured by the square root of the
$f$-divergence. For purposes of clarity, we write $M_{KL}(\epsilon,
S), M_C(\epsilon, S)$ and $M_l(\epsilon, S)$ for $M_f(\epsilon, S)$
when the function $f$ equals $x \log x$, $x^2 - 1$ and  $x^l - 1$ and respectively. 
\end{enumerate} 
We note here that the probability measures $Q_{\alpha}, \alpha \in G$
in the definition of $M_f(\epsilon, S)$ do not need to be included in
the set $\left\{P_{\theta}, \theta \in \Theta \right\}$ and the set
$G$ just denotes the index set and need not have any relation to $S$
or $\Theta$.  
\begin{theorem}\label{myb.thm}
The minimax risk $R$ satisfies the inequality $R \geq \sup_{\eta>0,
  \epsilon >0} \ell(\eta/2) (1-\star)$ where $\star$ stands for
any of the following quantities
\begin{equation}\label{myb}
 \frac{\log 2 + \log M_{KL}(\epsilon, \Theta) + \epsilon^2}{\log
    N(\eta)}
\end{equation}
\begin{equation}\label{myb.chi}
\frac{1}{N(\eta)} +
  \sqrt{\frac{(1+\epsilon^2)M_C(\epsilon, \Theta)}{N(\eta)}}
\end{equation}
and for $l > 1, l \neq 2$, 
\begin{equation}\label{myb.l}
\left( \frac{1}{N(\eta)^{l-1}} +
\frac{(1+\epsilon^2)M_l(\epsilon, \Theta)^{l-1}}{N(\eta)^{l-1}} \right)^{1/l}.
\end{equation}
\end{theorem}
In the sequel, by inequality~\eqref{myb.chi}, we mean the inequality
$R \geq \sup_{\eta >0, \epsilon >0} \ell(\eta/2)(1 - \star)$ with $\star$
representing~\eqref{myb.chi} and similarly for inequalities~\eqref{myb}
and~\eqref{myb.l}. 
\begin{IEEEproof}
We shall give the proof of inequality~\eqref{myb.chi}. The remaining
two inequalities are proved in a similar manner. 
Fix $\eta > 0$. By the definition of $N(\eta)$, one can find a
 finite subset $F \subset \Theta$ with cardinality $|F| \geq N(\eta)$
 such that $\rho(\theta, \theta ') \geq \eta$ for $\theta,  \theta
 ' \in F$ and $\theta \neq \theta '$. We then employ the inequality $R
 \geq \ell(\eta/2)  r$, where $r$ is defined as
 in~\eqref{rdef}. Inequality~\eqref{chifano} can now be used to obtain
\begin{equation*}
r  \geq  1 - \frac{1}{\sqrt{|F|}}
\sqrt{\frac{\inf_Q \sum_{\theta \in F}  \chi^2(P_{\theta}|| Q)}{|F|}} - \frac{1}{|F|}.
\end{equation*}
We now fix $\epsilon > 0$ and use the definition of $M_C(\epsilon,
F)$ to get a finite set $G$ with cardinality $\leq M_C(\epsilon, F)$
and probability measures $Q_{\alpha}, \alpha \in G$ such that
$\sup_{\theta \in S} \min_{\alpha \in G} \chi^2(P_{\theta}||Q_{\alpha})
\leq \epsilon^2$.
We then use inequality~\eqref{chijsd.eq} to get that
\begin{equation*}
\inf_Q \frac{1}{|F|} \sum_{\theta \in F} \chi^2(P_{\theta}||Q) \leq M_C(\epsilon, F)
\left( 1 + \epsilon^2 \right) - 1. 
\end{equation*}
The proof is complete by the trivial observation $M_C(\epsilon, F)
\leq M_C(\epsilon, \Theta)$. 
\end{IEEEproof}
The inequality~\eqref{myb} is due to Yang and Barron~\cite[Proof of
Theorem 1]{YangBarron}. In their paper, Yang and Barron mainly
considered the problem of estimation from $n$ independent and
identically distributed observations. However their method results in
inequality~\eqref{myb} which applies to every estimation
problem. Inequalities~\eqref{myb.chi} and~\eqref{myb.l} are new. 

Note that the lower bounds for $R$ given in
Theorem~\ref{myb.thm} all depend only on the 
quantities $N(\eta)$ and $M_f(\epsilon, \Theta)$, which describe
packing/covering properties of the entire parameter space
$\Theta$. Consequently, these inequalities only involve global metric entropy
properties. This is made possible by the use of
inequalities in Theorem~\ref{ybgen}. In applications of Fano's
inequality~\eqref{normfano} with the standard bound $J_f \leq
\max_{\theta, \theta ' \in F} D(P_{\theta}||P_{\theta '})$ as
well as in the application of other popular methods for obtaining
minimax lower bounds like Le Cam's method or Assouad's lemma, one
needs to construct the finite subset $F$ of the parameter space in a
very special way: the parameter values in $F$ should be reasonably
separated in the metric $\rho$ and also, the probability measures
$P_{\theta}, \theta \in F$ should be close in some probability
metric. In contrast, the application of Theorem~\ref{myb.thm}
does not require the construction of such a special subset $F$.

Yang and Barron~\cite{YangBarron} have successfully applied
inequality~\eqref{myb} to achieve minimax lower bounds of the optimal
rate for 
many nonparametric density estimation and regression problems where
$N(\eta)$ and $M_{KL}(\epsilon, \Theta)$ can be deduced from standard
results in approximation theory for function classes. We
refer the reader to~\cite{YangBarron} for
examples. In some of these examples, inequality~\eqref{myb.chi} can
also be applied to get optimal lower bounds. In section~\ref{suppest},
we shall employ inequality~\eqref{myb.chi} to obtain a new minimax
lower bound in the problem of reconstructing convex bodies from noisy
support function measurements.

But prior to that, let us assess the performance of
inequality~\eqref{myb.chi} in certain standard parametric estimation
problems. In these problems, an interesting contrast arises between
the two minimax lower bounds~\eqref{myb} and~\eqref{myb.chi}: the
inequality~\eqref{myb} only results in a sub-optimal lower bound on
the minimax risk (this observation, due to Yang and Barron~\cite[Page
1574]{YangBarron}, is also explained in Example~\ref{parnorm} below)
while~\eqref{myb.chi} produces rate-optimal lower bounds. 

Our intention here is to demonstrate, with the help of the subsequent three
examples, that inequality~\eqref{myb.chi} works even
for finite dimensional parametric estimation problems, a scenario in which
it is already known~\cite[Page 1574]{YangBarron} that inequality~\eqref{myb}
fails. Of course, obtaining optimal minimax rates in 
such problems is facile in most situations. For
example, a two-points argument based on Hellinger distance gives the
optimal rate, as is widely recognized since Le
Cam~\cite{LeCam:73AnnStat}. But the point here is that even in finite
dimensional situations, global metric entropy features are adequate
for obtaining rate-optimal minimax lower bounds. This is contrary to
the usual claim that in order to establish rate-optimal lower bounds
in parametric settings, one needs more information than global entropy
characteristics~\cite[Page 1574]{YangBarron}.  

In each of the ensuing three examples, we take the parameter space
$\Theta$ to be a bounded interval of the real line and we consider the problem of
estimating a parameter $\theta \in \Theta$ from $n$ independent
observations distributed according to $m_{\theta}$, where $m_{\theta}$
is a probability measure on the real line. The probability measure
$P_{\theta}$ accordingly equals the $n$-fold product of
$m_{\theta}$. We shall work with the squared error loss so that $\ell(x)
= x^2$, $\rho$ is the Euclidean distance on the real line and
$N(\eta)$ can be taken to $c_1/\eta$ for $\eta \leq \eta_0$ 
where $c_1$ and $\eta_0$ are positive constants depending on the
bounded parameter space alone. We shall encounter more positive
constants $c, c_2, c_3, c_4, c_5, \epsilon_0$ and $\epsilon_1$ in the
examples all of which depend possibly on the parameter space alone and
thus, independent of $n$. 
\begin{example}\normalfont
\label{parnorm}
Suppose that $m_{\theta}$ equals the normal distribution with mean $\theta$
and variance 1. It can be readily verified that, for $\theta, \theta '
\in \Theta$, one has
\begin{equation*}
  D(P_{\theta}||P_{\theta '}) = \frac{n}{2}|\theta - \theta '|^2
\end{equation*}
and
\begin{equation*}
 \chi^2(P_{\theta}||P_{\theta '}) = \exp \left(n|\theta
  - \theta '|^2\right)-1.
\end{equation*}
It follows directly that  $D(P_{\theta}||P_{\theta '}) \leq
\epsilon^2$ if and only if 
$|\theta - \theta '| \leq \sqrt{2} \epsilon/\sqrt{n}$ and
$\chi^2(P_{\theta}||P_{\theta '}) \leq \epsilon^2$ if and only if
$|\theta - \theta '| \leq \sqrt{\log(1+\epsilon^2)}/\sqrt{n}$. As a
result, we can take 
\begin{equation}\label{ucp}
M_{KL}(\epsilon, \Theta) = \frac{c_2
    \sqrt{n}}{\epsilon} \text{ and } M_C(\epsilon, \Theta) = \frac{c_2
  \sqrt{n}}{\sqrt{\log (1+\epsilon^2)}}
\end{equation}
for $\epsilon \leq \epsilon_0$. 
Now, inequality~\eqref{myb} says that the minimax risk $R_n$ satisfies 
\begin{equation*}
R_n  \geq  \sup_{\eta \leq \eta_0, \epsilon \leq \epsilon_0} \frac{\eta^2}{4}
\left( 1 - \frac{\log 2 + \log (c_2 \sqrt{n}/\epsilon) +
    \epsilon^2}{\log (c_1/\eta)}\right). 
\end{equation*}
The function $\epsilon \mapsto \epsilon^2 - \log \epsilon$ is
minimized on $[0, \epsilon_0]$ at, say, $\epsilon = \epsilon_1$ and we
then get 
\begin{equation}\label{noryb}
R_n  \geq  \sup_{\eta \leq \eta_0} ~\frac{\eta^2}{4} \left(1 -
  \frac{\log n + c_3}{2\log c_1 + 2 \log(1/\eta)} \right),
\end{equation}
where $c_3$ is a function of $c_2$ and $\epsilon_1$. We now note that
when $\eta = c/\sqrt{n}$ for a constant $c$, the 
quantity inside the parantheses on the right hand side
of~\eqref{noryb} converges to 0 as $n$ goes to $\infty$. This means
that inequality~\eqref{myb} only gives lower bounds of inferior order
for $R_n$, the optimal order being, of course, $1/n$.

On the other hand, we shall show below that 
inequality~\eqref{myb.chi} gives $R_n \geq c/n$ for a positive
constant $c$. Indeed, inequality~\eqref{myb.chi} says that 
\begin{equation*}
R_n \geq \sup_{\eta \leq \eta_0, \epsilon \leq \epsilon_0}~\frac{\eta^2}{4}
\left(1 - \frac{\eta}{c_1} - \sqrt{\eta
    \sqrt{n}}\sqrt{\frac{c_2(1+\epsilon^2)}{c_1\sqrt{\log(1+\epsilon^2)}}}\right).       
\end{equation*}
Taking $\epsilon = \epsilon_0$ and $\eta = c_3/\sqrt{n}$, we get
\begin{equation}
R_n  \geq  \frac{c_3^2}{4n} \left(1 -
  \frac{c_3}{c_1\sqrt{n}} - c_4 \sqrt{c_3} \right),
\end{equation}
where $c_4$ depends only on $c_1, c_2$ and $\epsilon_0$. Hence by
choosing $c_3$ small, we get that $R_n \geq c/n$ for all large $n$.  
\end{example}
\begin{example}\normalfont
Suppose that $\Theta$ is a compact interval of the positive real line
that is bounded away from 
zero and suppose that $m_{\theta}$ denotes the uniform distribution on
$[0, \theta]$. It is then elementary to check that the chi-squared
divergence between $P_{\theta}$ and $P_{\theta '}$ equals $(\theta
'/\theta)^n-1$ if $\theta \leq \theta '$ and $\infty$ otherwise. It
follows accordingly that $\chi^2(P_{\theta}||P_{\theta '}) \leq
\epsilon^2$ provided 
\begin{equation}\label{mj}
  0 \leq \theta ' - \theta \leq \frac{\theta \log (1+\epsilon^2)}{n}. 
\end{equation}
Because the parameter space is a compact interval bounded away from
zero, in order to ensure~\eqref{mj}, it is enough to require that $0
\leq \theta ' - \theta \leq c_2 \log(1+\epsilon^2)/n$. Therefore, we
can take
\begin{equation*}
  M_C(\epsilon, \Theta) = \frac{c_3n}{\log(1+\epsilon^2)}
\end{equation*}
for $\epsilon \leq \epsilon_0$. Inequality~\eqref{myb.chi} now implies
that 
\begin{equation*}
  R_n \geq \sup_{\eta \leq \eta_0, \epsilon \leq \epsilon_0} \frac{\eta^2}{4} \left(1 -
    \frac{\eta}{c_1} - \sqrt{\eta n}
    \sqrt{\frac{c_3(1+\epsilon^2)}{c_1\log(1+\epsilon^2)}}  \right). 
\end{equation*}
Taking $\epsilon = \epsilon_0$ and $\eta = c_4/n$, we get that
\begin{equation*}
  R_n \geq \frac{c_4^2}{4n^2} \left(1 - \frac{c_4}{nc_1} - \sqrt{c_4}
    c_5 \right),
\end{equation*}
where $c_5$ depends only on $c_1, c_3$ and $\epsilon_0$. 
Hence by choosing $c_4$ sufficiently small, we get that $R_n \geq
c/n^2$ for all large $n$. This is the optimal minimax rate for this
problem as can be seen by estimating $\theta$ by the maximum of the
observations. 
\end{example}
\begin{example}\normalfont
Suppose that $m_{\theta}$ denotes the uniform distribution on the
interval $[\theta, \theta + 1]$. We shall argue that $M_C(\epsilon,
\Theta)$ can be chosen to be 
  \begin{equation}\label{lara}
    M_C(\epsilon, \Theta) = \frac{c_2}{(1+\epsilon^2)^{1/n}-1}
  \end{equation}
for a positive constant $c_2$ at least for large $n$. To see this, let
us define $\epsilon '$ so that $2\epsilon
':=(1+\epsilon^2)^{1/n}-1$ and let $G$ denote an $\epsilon '$-grid of 
points in the interval $\Theta$; $G$ would contain at most
$c_2/\epsilon '$  points when $\epsilon \leq \epsilon_0$. For a point
$\alpha$ in the grid, let $Q_{\alpha}$ denote the $n$-fold product of
the uniform distribution on the interval $[\alpha, \alpha + 1 +
2\epsilon ']$. Now, for a fixed $\theta \in \Theta$, let $\alpha$
denote the point in the grid such that $\alpha \leq \theta \leq \alpha
+ \epsilon '$. It can then be checked that the chi-squared divergence
between $P_{\theta}$ and $Q_{\alpha}$ is equal to $(1+2 \epsilon ')^n
- 1 = \epsilon^2$. Hence $M_C(\epsilon, \Theta)$ can be taken to be
the number of probability measures $Q_{\alpha}$, which is the same as the
number of points in $G$. We thus have~\eqref{lara}. It can be checked
by elementary calculus (Taylor expansion, for example) that the
inequality  
\begin{equation*}
  (1+\epsilon^2)^{1/n} - 1 \geq \frac{\epsilon^2}{n} -
  \frac{1}{2n}\left( 1 - \frac{1}{n}\right)\epsilon^4
\end{equation*}
holds for $\epsilon \leq \sqrt{2}$ (in fact for all $\epsilon$, but
for $\epsilon > \sqrt{2}$, the right hand side above may be
negative). Therefore for $\epsilon \leq \min(\epsilon_0, \sqrt{2})$,
we get that 
\begin{equation*}
  M_C(\epsilon, \Theta) \leq \frac{2nc_2}{2\epsilon^2 - (1-1/n)\epsilon^4}.
\end{equation*}
From inequality~\eqref{myb.chi}, we get that for every $\eta \leq
\eta_0$ and $\epsilon \leq \min(\epsilon_0, \sqrt{2})$, 
\begin{equation*}
  R_n \geq \frac{\eta^2}{4} \left(1 - \frac{\eta}{c_1} - \sqrt{n \eta}
  \sqrt{\frac{2(1+\epsilon^2)c_2}{c_1 \left(2\epsilon^2 -
        (1-1/n)\epsilon^4 \right)}}\right).
\end{equation*}
If we now take $\epsilon = \min(\epsilon_0, 1)$ and $\eta = c_3/n$, we
see that the quantity inside the parantheses converges to
$1-\sqrt{c_3}c_4$ where $c_4$ depends only on $c_1, c_2$ and
$\epsilon_0$. Therefore by choosing $c_3$ sufficiently small, we get
that $R_n \geq c/n^2$. This is the optimal minimax rate for this
problem as can be seen by estimating $\theta$ by the minimum of the
observations. 
\end{example}
The fact that inequality~\eqref{myb.chi} produced optimal lower bounds for
the minimax risk in each of the above three examples is reassuring but
not really exciting because, as we mentioned before, there are other simpler
methods of obtaining such bounds in these examples. We presented them
as simple toy examples to evaluate the performance of~\eqref{myb.chi},
to present a difference between~\eqref{myb} and~\eqref{myb.chi} (which
provides a justification for using divergences other than the
Kullback-Leibler divergence for lower bounds) and also to stress the
fact that global packing and covering characteristics 
are enough to obtain optimal minimax lower bounds. In order to
convince the reader of the effectiveness of~\eqref{myb.chi} in more
involved situations, we now apply it to obtain the optimal minimax
rate in a $d$-dimensional normal mean estimation problem. We
are grateful to an anonymous referee for communicating this example to
us. Another non-trivial application of~\eqref{myb.chi} is presented in
the next section.  
\begin{example}\normalfont
\label{ddimnorm}
Let $\Theta$ denote the ball in $\R^d$ of radius
$\Gamma$ centered at the origin. Let us consider the problem of estimating
$\theta \in \Theta$ from an observation $X$ distributed according to
the normal distribution with mean $\theta$ and variance covariance
matrix $\sigma^2 I_d$, where $I_d$ denotes the 
identity matrix of order $d$. Thus $P_{\theta}$ denotes the $N(\theta,
\sigma^2 I_d)$ distribution. We assume squared error loss so that 
$\ell(x) = x^2$ and $\rho$ is the Euclidean distance on $\R^d$. 

We shall use inequality~\eqref{myb.chi} to show that the minimax risk
$R$ for this problem is larger than or equal to a constant multiple of
$d \sigma^2$ when $\Gamma \geq \sigma \sqrt{d}$. We begin by arguing that
we can take
\begin{equation}\label{pete}
  N(\eta) = \left( \frac{\Gamma}{\eta}\right)^d,
  M_C(\epsilon, \Theta) = \left(\frac{3\Gamma}{\sigma
      \sqrt{\log(1+\epsilon^2)}} \right)^d
\end{equation}
whenever $\sigma \sqrt{\log(1+\epsilon^2)} \leq \Gamma$. 

For $N(\eta)$, we first note that the $\eta$-packing number
of the metric space $(\Theta, \rho)$ is bounded from below by its
$\eta$-covering number. Now, for any $\eta$-covering set, the space
$\Theta$ is contained in the union of the balls of radius $\eta$ with
centers in the covering set and hence the volume of $\Theta$ must be
smaller than the sum of the volumes of these balls. Therefore, the
number of points in the $\eta$-covering set must be at least
$(\Gamma/\eta)^d$. Since this is true for every $\eta$-covering set,
it follows that the $\eta$-covering number and hence the
$\eta$-packing number is not smaller than $(\Gamma/\eta)^d$.

For $M_C(\epsilon, \Theta)$, we first observe that for $\theta, \theta
' \in \Theta$, the chi-squared divergence between 
$P_{\theta}$ and $P_{\theta '}$ can be easily computed (because they
are normal distributions with the same covariance matrix) to be
$\chi^2(P_{\theta}||P_{\theta '}) = \exp \left(\rho^2(\theta, \theta
  ')/\sigma^2 \right) - 1$. Therefore $\chi^2(P_{\theta}||P_{\theta
  '}) \leq \epsilon^2$ if and only if $\rho(\theta, \theta ') \leq
\epsilon ' := \sigma \sqrt{\log(1+\epsilon^2)}$. As a result,
$M_C(\epsilon, \Theta)$ can be taken to be any upper bound on the
$\epsilon '$-covering number of $(\Theta, \rho)$. The $\epsilon
'$-covering number, as noted previously, is bounded from above by the
$\epsilon '$-packing number. Now, for any $\epsilon '$-packing set,
the balls of radius $\epsilon '/2$ with centers in the packing set are
all disjoint and their union is contained in the ball of radius
$\Gamma + (\epsilon '/2)$ centered at the origin. Consequently, the
sum of the volumes of these balls is smaller than the volume of the
ball of radius $\Gamma + (\epsilon '/2)$ centered at the
origin. Therefore, the number of points in the $\epsilon '$-packing
set is at most $(1+(2\Gamma/\epsilon '))^d \leq (3\Gamma/\epsilon
')^d$ provided $\epsilon ' \leq \Gamma$. Since this is true for every
$\epsilon '$-packing set, it follows that the $\epsilon '$-packing
number and hence the $\epsilon '$-covering number is not larger than
$(3\Gamma/\epsilon ')^d$.

We can thus apply inequality~\eqref{myb.chi} with~\eqref{pete} to get
that, for every $\eta > 0$ and $\epsilon > 0$ such that $\sigma
\sqrt{\log (1+\epsilon^2)} \leq \Gamma$, we have
\begin{equation*}
  R \geq \frac{\eta^2}{4} \left(1 - \left(\frac{\eta}{\Gamma}\right)^d
    - \left(\frac{3\eta}{\sigma}\right)^{d/2}
    \frac{\sqrt{1+\epsilon^2}}{(\log (1+\epsilon^2))^{d/4}} \right). 
\end{equation*}
Now by elementary calculus, it can be checked that the function
$\epsilon \mapsto \sqrt{1+\epsilon^2}/(\log (1+\epsilon^2))^{d/4}$ is
minimized (subject to $\sigma \sqrt{\log(1+\epsilon^2)} \leq \Gamma$)
when $1+\epsilon^2 = e^{d/2}$. We then get that
\begin{equation*}
  R \geq \sup_{\eta > 0} \frac{\eta^2}{4} \left(1 -
    \left(\frac{\eta}{\Gamma}\right)^d -
    \left(\frac{18e\eta^2}{\sigma^2d} \right)^{d/4}\right).  
\end{equation*}
We now take $\eta = c_1 \sigma \sqrt{d}$ and since $\Gamma \geq \sigma
\sqrt{d}$, we obtain
\begin{equation*}
  R \geq \frac{c_1^2\sigma^2d}{4} \left(1 - c_1^d - (18ec_1^2)^{d/4} \right).
\end{equation*}
We can therefore choose $c_1$ small enough (independent of $d$) to
obtain that $R \geq cd\sigma^2$. Note that, up to constants, this
lower bound is optimal for $R$ because $\E \rho^2(X, \theta) =
d\sigma^2$. 
\end{example}
\section{Reconstruction of convex bodies from noisy  support function
  measurements}\label{suppest} 
In this section, we shall present a novel application of the global minimax
lower bound~\eqref{myb.chi}. Let $d \geq 2$ and let $K$ be a convex
body in $\R^d$, i.e., $K$ is compact, convex and has a non-empty
interior. The support function of $K$, $h_K : S^{d-1} \rightarrow \R$,
is defined by 
\begin{equation*}
h_K(u)  :=  \sup \left\{\left<x, u \right> : x \in K
\right\} \text{ for } u \in S^{d-1}, 
\end{equation*}
where $S^{d-1} := \left\{ x \in \R^d : \sum_i x_i^2 = 1
\right\}$ is the unit sphere. We direct the reader to~\cite[Section
1.7]{Schneider} or~\cite[Section
13]{Rockafellar70book} for basic properties 
of support functions. An important property is that the support
function uniquely determines the convex body, i.e., $h_K = h_L$
if and only if $K = L$.

Let $\left\{u_i, i \geq 1 \right\}$ be a sequence of $d$-dimensional
unit vectors. Gardner, Kiderlen and
Milanfar~\cite{GardnerKiderlenMilanfar} (see their paper for earlier
references) considered the problem of reconstructing an unknown convex
body $K$ from noisy measurements of $h_K$ in the directions $u_1,
\dots, u_n$. More precisely, their problem was to estimate $K$ from
observations $Y_1, \dots, Y_n$ drawn according to the model $Y_i= 
h_K(u_i) + \xi_i, i = 1, \dots, n$ where $\xi_1, \dots, \xi_n$ are
independent and identically distributed mean zero gaussian random
variables. They constructed a convex body (estimator) $\hat{K}_n =
\hat{K}_n(Y_1, \dots, Y_n)$ having the property that, for
\textit{nice} sequences $\left\{ u_i, i \geq 1 \right\}$, the $L^2$
norm $||h_K - h_{\hat{K}_n}||_2$ (see~\eqref{lpnorm} below) converges
to zero at the rate $n^{-2/(d+3)}$ for dimensions $d = 2, 3, 4$ and at
a slower rate for dimensions $d \geq 5$ (see~\cite[Theorem
6.2]{GardnerKiderlenMilanfar}).  

We shall show here that in the same setting, it is impossible in the
minimax sense to construct estimators 
for $K$ converging at a rate faster than $n^{-2/(d+3)}$. This
implies that the least squares estimator
in~\cite{GardnerKiderlenMilanfar} is rate optimal for 
dimensions $d = 2, 3, 4$.  We shall need some
notation to describe our result. 

Let $\K^d$ denote the set of all convex bodies in $\R^d$ and
for $\Gamma > 0$, let $\K^d(\Gamma)$ denote the set of all convex bodies
in $\R^d$ that are contained in the closed ball of radius $\Gamma$ centered
at the origin so that $\K^d(1)$ denotes the set of all convex bodies
contained in the unit ball, which we shall denote by $B$. Note that
estimating $K$ is equivalent to estimating 
the function $h_K$ because the support function uniquely determines
the convex body. Thus we shall focus on the problem of estimating
$h_K$.

An estimator for $h_K$ is allowed to be a bounded
function on $S^{d-1}$ that depends on the data $Y_1, \dots, Y_n$. The
loss functions that we shall use are the $L^p$ norms for $p \in [1,
\infty]$ defined by
\begin{equation}\label{lpnorm}
||h_K - \hat{h}||_p := \left(\int_{S^{d-1}} |h_K(u) -
    \hat{h}(u)|^p du
  \right)^{1/p}
\end{equation}
for $p \in [1, \infty)$ and $||h_K - \hat{h}||_{\infty} := \sup_{u \in
  S^{d-1}} |h_K(u) - \hat{h}(u)|$. For convex bodies $K$ and $L$ and
$p \in [1, \infty]$, we shall also write $\delta_p(K, L)$ for $||h_K -
h_L||_p$ and refer to $\delta_p$ as the $L^p$ distance between $K$ and
$L$.

We shall consider the minimax risk of the problem of estimating $h_K$ from 
$Y_1, \dots, Y_n$ when $K$ is assumed to belong to $\K^d(\Gamma)$ i.e., we
are interested in the quantity  
\begin{equation*}
r_n(p, \Gamma)   :=  \inf_{\hat{h}} \sup_{K \in \K^d(\Gamma)}  \E_K||h_K -
\hat{h}(Y_1, \dots, Y_n)||_p.
\end{equation*}
The following is the main theorem of this section. 
\begin{theorem}\label{suppthm}
Fix $p \in [1, \infty)$ and $\Gamma > 0$. Suppose the errors $\xi_1, \dots,
\xi_n$ are independent normal random variables with mean zero and
variance $\sigma^2$. Then the minimax risk $r_n(p, \Gamma)$ satisfies
\begin{equation}\label{suppthm.eq}
r_n(p, \Gamma) \geq c  \sigma^{4/(d+3)} \Gamma^{(d-1)/(d+3)} n^{-2/(d+3)}, 
\end{equation}
for a constant $c$ that is independent of $n$. 
\end{theorem}
\begin{remark}\normalfont
In the case when $p = 2$, Gardner, Kiderlen and
Milanfar~\cite{GardnerKiderlenMilanfar} showed 
that the least squares estimator converges at the rate given by the
right hand side of~\eqref{suppthm.eq} for dimensions $d = 2, 3,
4$. Thus, at least for $p=2$, the lower bound given
by~\eqref{suppthm.eq} is optimal for dimensions $d = 2, 3, 4$. 
\end{remark}
We shall use inequality~\eqref{myb.chi} to
prove~\eqref{suppthm.eq}. First, let us put the support function 
estimation problem in the general estimation setting of the last
section. Let $\Theta := \left\{h_K : K \in \K^d(\Gamma) \right\}$ and let
$\ac$ be the collection of all bounded functions on the unit sphere
$S^{d-1}$. The metric $\rho$ on $\ac$ is just the $L^p$ norm and
$\ell(x) = x$. 

Finally, let $\samp = \R^n$ and for $f \in \Theta$, let $P_f$ be the
$n$-variate normal distribution with mean vector $(f(u_1), \dots,
f(u_n))$ and variance-covariance matrix $\sigma^2 I_n$, where $I_n$ is
the identity matrix of order $n$. 

In order to apply inequality~\eqref{myb.chi}, we need to determine
$N(\eta)$ and $M_C(\epsilon, \Theta)$. The quantity $N(\eta)$ is a
lower bound on the $\eta$-packing number of the set $\K^d(\Gamma)$ under
the $L^p$ norm. When $p = \infty$, Bronshtein~\cite[Theorem 4 and
Remark 1]{Bronshtein76} proved that there exist positive constants
$c'$ and $\eta_0$ depending only on $d$ such that $\exp\left( c'
  (\eta/\Gamma)^{(1-d)/2}\right)$ is a lower bound for the $\eta$-packing
number of $\Theta$ for $\eta \leq \eta_0$. It is a standard fact that
$p = \infty$ corresponds to the Hausdorff metric on $\K^d(\Gamma)$. 

It turns out that Bronshtein's result is actually true for every $p
\in [1, \infty]$ and not just for $p = \infty$. However, to the best
of our knowledge, this has not been proved anywhere in the literature. By
modifying Bronshtein's proof appropriately and using the
Varshamov-Gilbert lemma (see for example~\cite[Lemma
4.7]{Massart03Flour}), we provide, in Theorem~\ref{Rpackresult}, a
proof of this fact. Therefore from Theorem~\ref{Rpackresult}, we can
take   
\begin{equation}\label{Neta}
\log N(\eta)  = c' \left(\frac{\Gamma}{\eta} \right)^{(d-1)/2}
\text{ for } \eta \leq \eta_0, 
\end{equation} 
where $c'$ and $\eta_0$ are constants depending only on $d$ and $p$.

Now let us turn to $M_C(\epsilon, \Theta)$. For $f, g \in \Theta$,
$P_f$ and $P_g$ are normal distributions with the same covariance
matrix and hence the chi-squared divergence between $P_f$ and $P_g$
can be seen to be
\begin{align*}
\chi^2(P_f||P_g) & = \exp \left[ \frac{1}{\sigma^2} \sum_{i=1}^n \left(f(u_i) - g(u_i)
  \right)^2 \right] - 1 \\
&\leq  \exp \left[ \frac{n||f-g||^2_{\infty}}{\sigma^2} \right] - 1.
\end{align*}
It follows that 
\begin{equation}\label{klhaus}
||f - g||_{\infty} \leq \epsilon '\Longrightarrow  \chi^2(P_f||P_g) \leq \epsilon^2.
\end{equation}
where $\epsilon ':= \sigma \sqrt{\log (1 +
  \epsilon^2)}/\sqrt{n}$. Let $W_{\epsilon '}$ be the smallest $W$ for
which there exist sets $K_1, \dots, K_W$ in $\K^d(\Gamma)$ having the
property that for every set $K 
\in \K^d(\Gamma)$, there exists a $K_j$ such that the Hausdorff distance
between $K$ and $K_j$ is less than or equal to $\epsilon '$. It must
be clear from~\eqref{klhaus} that $M_C(\epsilon, \Theta)$ can be taken to be
a number larger than $W_{\epsilon '}$. Bronshtein~\cite[Theorem
3 and Remark 1]{Bronshtein76} showed that there exist positive
constants $c''$ and $\epsilon_0$ depending only on $d$ such that  
\begin{equation*}
\log W_{\epsilon '}  \leq c'' \left(\frac{\Gamma}{\epsilon '}
\right)^{(d-1)/2} \text{ for } \epsilon ' \leq \epsilon_0. 
\end{equation*}
Hence for all $\epsilon$ such that $\log(1 + \epsilon^2) \leq n \epsilon_0^2/\sigma^2$, we
can take
\begin{equation}\label{Meps}
\log M_C(\epsilon, \Theta)  = c'' \left(\frac{\Gamma \sqrt{n}}{\sigma 
    \sqrt{\log(1 + \epsilon^2)}} \right)^{(d-1)/2}.
\end{equation}
We are now ready to prove inequality~\eqref{suppthm.eq}. We shall
define two quantities 
\begin{equation*}
\eta(n):= c \sigma^{4/(d+3)} \Gamma^{(d-1)/(d+3)} n^{-2/(d+3)} 
\end{equation*}
and
\begin{equation*}
u(n) := \left(\frac{\Gamma\sqrt{n}}{\sigma}
\right)^{(d-1)/(d+3)}. 
\end{equation*}
where $c = c(d, p)$ will be specified shortly. Also let
$\epsilon(n)$ be such that $\log (1 + \epsilon^2(n)) = u^2(n)$. 
Clearly as $n \rightarrow \infty$, we have $\eta(n) \rightarrow 0$,
$u(n) \rightarrow \infty$ and $u(n)/\sqrt{n} \rightarrow
0$. As a result $\eta(n) \leq \eta_0$ and
$u^2(n) \leq n\epsilon_0^2/\sigma^2$ for large $n$ and therefore
from~\eqref{Neta} and~\eqref{Meps}, we get that
\begin{equation*}
  \log N(\eta(n))  = c' \left(\frac{\Gamma}{\eta(n)} \right)^{(d-1)/2}
  = \frac{c'}{c^{(d-1)/2}} u^2(n).
\end{equation*}
and
\begin{equation*}
  \log M_C(\epsilon(n), \Theta)  = c'' \left(\frac{\Gamma
      \sqrt{n}}{\sigma u(n)} \right)^{(d-1)/2} = c'' u^2(n).
\end{equation*}
We now apply inequality~\eqref{myb.chi} (recall that $\ell(x) = x$) to
obtain that $r_n(p, \Gamma)$ is bounded from below by  
\begin{equation*}
  \frac{\eta(n)}{2} \left[1 - \frac{1}{N(\eta(n))} - \exp
    \left(\frac{u^2(n)}{2} \left(1 + c'' - \frac{c'}{c^{(d-1)/2}}
      \right) \right) \right]
\end{equation*}
for all large $n$. If we now choose $c$ so that $c^{(d-1)/2} =
c'/(2+2c'')$, we get that
\begin{equation*}
  r_n(p, \Gamma) \geq \frac{\eta(n)}{2} \left[1 - \frac{1}{N(\eta(n))}
  - \exp \left(\frac{-u^2(n)}{2}(1+c'') \right)\right].
\end{equation*}
Now observe that as $n \rightarrow \infty$, the quantity $\eta(n)$
goes to 0 and hence $N(\eta(n))$ goes to $\infty$. Further, as we have
already noted, $u(n)$ goes to $\infty$. It follows hence that $r_n(p,
\Gamma) \geq \eta(n)/4$ for all large $n$. By choosing $c$ even smaller, we
can make inequality~\eqref{suppthm.eq} true for all $n$. 
\section{A covariance matrix estimation example}\label{covmat}
In the previous section, we have used the global minimax lower
bound~\eqref{myb.chi}. However, in some situations, the global entropy
numbers might be difficult to bound. In such cases,
inequalities~\eqref{myb} and~\eqref{myb.chi} are, of course, not
applicable and we are unaware of the use of inequality~\eqref{jsyb} in
conjuction with Fano's inequality~\eqref{normfano} in the
literature. The standard examples use~\eqref{normfano} with the bound
$J_f \leq \min_{\theta, \theta ' \in F}D(P_{\theta}||P_{\theta '})$
while the examples in~\cite{YangBarron} all deal with the
case when global entropies are available. In this section, we shall
demonstrate how a recent minimax lower bound due to Cai, Zhang and
Zhou~\cite{CaiZhangZhou2009} can also be proved using
inequalities~\eqref{normfano} and~\eqref{jsyb}.  

Cai, Zhang and Zhou~\cite{CaiZhangZhou2009} considered $n$ independent $p
\times 1$ random vectors $X_1, \dots, X_n$ distributed according to
$N_p(0, \Sigma)$. Suppose that the entries of the $p \times p$
covariance matrix $\Sigma = (\sigma_{ij})$ decay at a certain rate as
we move away from the diagonal. Specifically, let us suppose that for
a fixed positive constant $\alpha > 0$, the entries $\sigma_{ij}$ of
$\Sigma$ satisfy the inequality $\sigma_{ij} \leq  |i-j|^{-\alpha-1}$ for $i
\neq j$. Cai, Zhang and Zhou~\cite{CaiZhangZhou2009} showed that when $p$ is large compared
to $n$, it is impossible to estimate $\Sigma$ from $X_1, \dots, X_n$
in the spectral norm at a rate faster than $n^{-\alpha/(2\alpha +
  1)}$. More precisely, they showed that when $p \geq Cn^{1/(2\alpha +
  1)}$, 
\begin{equation}
  \label{covres}
R_n(\alpha)  :=   \inf_{\hat{\Sigma}} \sup_{\Sigma \in \Theta}
\E_{\Sigma}||\hat{\Sigma} -
  \Sigma ||  \geq  c~n^{-\alpha/(2\alpha + 1)},
\end{equation}
where $c$ and $C$ denote positive constants depending only on $\alpha$. Here
$\Theta$ denotes the collection of all covariance matrices
$\Sigma = (\sigma_{ij})$ satisfying  $\sigma_{ij} \leq  |i-j|^{-\alpha-1}$ for $i
\neq j$ and the norm $||.||$ is the spectral norm (largest
eigenvalue). 

Cai, Zhang and Zhou~\cite{CaiZhangZhou2009} used Assouad's lemma
for the proof of the inequality~\eqref{covres}. We shall use
inequalities~\eqref{normfano} and~\eqref{jsyb}. Moreover, the choice of
the finite subset $F$ that we use is different from the one used
in~\cite[Equation (17)]{CaiZhangZhou2009}. This 
makes our approach different from the general method, due to
Yu~\cite{Yu97lecam}, of replacing Assouad's lemma by Fano's inequality.

Throughout, $\Delta$ denotes a constant that depends on $\alpha$
alone. The value of the constant might vary from place to place.

Consider the matrix $A = (a_{ij})$ with $a_{ij} = 1$ for $i=j$ 
and $a_{ij} = 1/(\Delta|i-j|^{\alpha+1})$ for $i \neq j$. For $\Delta$
sufficiently large (depending on $\alpha$ alone), $A$ is positive
definite and belongs to $\Theta$. Let 
us fix a positive integer $k \leq p/2$ and partition $A$ as
\begin{equation*}
A  =  \left[\begin{array}{c|c} A_{11} & A_{12}
     \\ \hline \rule[13pt]{0pt}{0pt} A_{12}^T &
    A_{22} \end{array} \right],
\end{equation*}
where $A_{11}$ is $k \times k$ and $A_{22}$ is $(p-k) \times (p-k)$.
For each $\tau \in \R^k$, we define the matrix
\begin{equation*}
  A(\tau)  := \left[\begin{array}{c|c} A_{11} & A_{12}(\tau)
     \\ \hline \rule[13pt]{0pt}{0pt} \left(A_{12}(\tau)\right)^T &
    A_{22} \end{array} \right],
\end{equation*}
where $A_{12}(\tau)$ is the $k \times (p-k)$ matrix obtained by
premultiplying $A_{12}$ with the $k \times k$ diagonal matrix with
diagonal entries $\tau_1, \dots, \tau_k$. Clearly, $A(\tau) \in
\Theta$ for all $\tau \in \left\{0, 1 \right\}^k$. We shall need the
following two lemmas in order to prove inequality~\eqref{covres}. 
\begin{lemma}\label{hamlb}
  For $\tau, \tau '\in \left\{0, 1 \right\}^k, \tau \neq \tau '$, we have
\begin{equation}
  \label{hamlb.eq}
  ||A(\tau) - A(\tau ')||  \geq 
 \frac{1}{\Delta k^{\alpha}} \sqrt{\frac{\Upsilon(\tau, \tau ')}{k}},  
\end{equation}
where $\Upsilon(\tau, \tau ') := \sum_{r=1}^k \left\{\tau_r \neq
  \tau'_r \right\}$ denotes the Hamming distance between $\tau$ and
$\tau'$.    
\end{lemma}
\begin{IEEEproof}
Fix $\tau, \tau' \in \left\{0, 1 \right\}^k$ with $\tau \neq \tau'$.
Let $v$ denote the $p \times 1$ vector 
$\left(0_k, 1_k, 0_{p-2k} \right)^T$, where $0_k$ denotes the
$k\times 1$ vector of zeros etc. Clearly $||v||^2 = k$ and $(A(\tau) -
A(\tau ')) v$ will be a vector of the form $(u, 0)^T$ for
some $k \times 1$ vector $u = (u_1, \dots, u_k)^T$. Moreover $u_r =
\sum_{s=1}^k (\tau_r - \tau_r') a_{r, k+s}$ and hence 
\begin{align*}
  |u_r|  &= \frac{ \left\{ \tau_r \neq \tau_r' \right\}}{\Delta}
  \sum_{s=1}^k \frac{1}{|r-k-s|^{\alpha+1}} \\
& \geq \frac{\left\{ \tau_r
    \neq \tau_r' \right\}}{\Delta} \sum_{i=k}^{2k -1} \frac{1}{i^{\alpha+1}}
  \geq \frac{\left\{ \tau_r \neq \tau_r'
    \right\}}{\Delta}\frac{1}{k^{\alpha}}. 
\end{align*}
Therefore, 
\begin{equation*}
  ||\left(A(\tau) - A(\tau ')\right)v||^2 \geq \sum_{r=1}^k u_r^2 \geq
  \frac{1}{\Delta^2 k^{2 \alpha}} \Upsilon(\tau, \tau '). 
\end{equation*}
The proof is complete because $||v||^2 = k$. 
\end{IEEEproof}
\begin{lemma}\label{covublem}
  Let $1 \leq m < k, \tau \in \left\{0, 1 \right\}^k$ and  $\tau'
  := (0, \dots, 0, \tau_m, \dots, \tau_k)$. Then
  \begin{equation*}
    D \left(N(0, A(\tau)) || N(0, A(\tau')) \right)  \leq
    \frac{\Delta}{(k-m)^{2 \alpha}}. 
  \end{equation*}
\end{lemma}
\begin{IEEEproof}
The key is to note that one has the inequality $D \left(N(0, A(\tau))
  || N(0, A(\tau')) \right) 
  \leq \Delta ||A(\tau) - A(\tau')||^2_F$, where $||A||_F := \left(\sum_{i,
    j} a_{ij}^2 \right)^{1/2}$ denotes the Frobenius norm. The proof
of this assertion can be found in~\cite[Proof of
Lemma 6]{CaiZhangZhou2009}. We can now bound
\begin{align*}
  ||A(\tau) - A(\tau')||^2_F &\leq 2\sum_{r=1}^{m-1} \tau_r^2
  \sum_{j=1}^{p-k} a_{r, k+j}^2 \\
&\leq \Delta \sum_{r=1}^{m-1}\sum_{j=1}^{p-k} \frac{1}{|r-k-j|^{2\alpha+2}}
\\
&\leq \Delta  \sum_{r=1}^{m-1} \sum_{j=1}^{\infty}
\frac{1}{|k-r+j|^{2\alpha+2}}\\
& \leq \Delta \sum_{r=1}^{m-1}
\frac{1}{(k-r)^{2\alpha + 1}} \leq \frac{\Delta}{(k-m)^{2\alpha}}.
\end{align*}
The proof is complete. 
\end{IEEEproof}
The Varshamov-Gilbert lemma (see for example~\cite[Lemma
4.7]{Massart03Flour}) asserts the existence of a subset
$W$ of $\left\{0, 1\right\}^k$ with $|W| \geq \exp(k/8)$ such that
$\Upsilon(\tau, \tau') \geq k/4$ for all $\tau, \tau' \in W$ with $\tau
\neq \tau'$. Let $F := \left\{A(\tau) : \tau \in W \right\}$. From
inequality~\eqref{normfano} and Lemma~\ref{hamlb}, we get that
\begin{equation}\label{emp2}
  R_n(\alpha)  \geq  \frac{1}{\Delta} \frac{1}{k^{\alpha}}
  \left(1 - \frac{\log 2 + \frac{1}{|W|}\sum_{A \in
        F}D(P_A||\bar{P})}{k/8} \right),
\end{equation}
where $P_A$ denotes the $n$-fold product of the $N(0, A)$ probability
measure and $\bar{P} := 
\sum_{A\in F} P_A/|W|$. Now for $1 \leq m < k$ and for $t \in
\left\{0, 1 \right\}^{k-m+1}$, let $Q_t$ denote the $n$-fold product
of the $N(0, A(0, \dots, 0, t_1, \dots, t_{k-m+1}))$ probability
measure. Applying inequality~\eqref{jsyb}, we get the quantity
$\sum_{A \in F} D(P_A||\bar{P})/|W|$ is bounded from above by
\begin{equation*}
 (k-m+1) \log 2 +
  \max_{A \in F} \min_{t \in \left\{ 0, 1\right\}^{k-m+1}} D(P_A||Q_t).
\end{equation*}
Now we use Lemma~\ref{covublem} to obtain
\begin{equation*}
\frac{1}{|W|}  \sum_{A \in F}D(P_A||\bar{P}) 
 \leq   \Delta \left[(k-m) + \frac{n}{(k-m)^{2\alpha}} \right]. 
\end{equation*}
Using the above in~\eqref{emp2}, we get
\begin{equation*}
  R_n(\alpha)  \geq  \frac{1}{\Delta} \frac{1}{k^{\alpha}}
  \left[1 - \frac{\Delta}{k}\left((k-m) + \frac{n}{(k-m)^{\alpha}}
    \right) \right]. 
\end{equation*}
Note that the above lower bound for $R_n(\alpha)$ depends on $k$ and
$m$, which are constrained to satisfy $2k \leq p$ and $1 \leq m <
k$. To get the best lower bound, we need to optimize the right hand
side of the above inequality over $k$ and $m$. It should be obvious
that in order to prove~\eqref{covres}, it is enough to take $k-m =
n^{1/(2 \alpha + 1)}$ and $k = 4\Delta n^{1/(2\alpha + 1)}$. The condition
$2k \leq p$ will be satisfied if 
$p \geq Cn^{1/(2\alpha + 1)}$ for a large enough $C$. It is elementary
to check that with these choices of $k$ and $m$,
inequality~\eqref{covres} is established.  
\section{A Packing Number Lower Bound}
In this section, we shall prove that for every $p \in [1, \infty]$ the
$\eta$-packing number $N(\eta; p, \Gamma)$ of $\K^d(\Gamma)$ under the $L^p$
metric is at least $\exp \left(c (\eta/\Gamma)^{(1-d)/2} \right)$ for a positive
$c$ and sufficiently small $\eta$. 
This means that there exist at least $\exp \left(c (\eta/\Gamma)^{(1-d)/2}
\right)$ sets in $\K^d(\Gamma)$ separated by at least $\eta$ in the $L^p$
metric. This result was needed in the proof of
Theorem~\ref{suppthm}. Bronshtein~\cite[Theorem 4 and Remark
1]{Bronshtein76} proved this for $p = \infty$ (the case of the
Hausdorff metric).  
\begin{theorem}\label{Rpackresult}
Fix $p \in [1, \infty]$. There exist positive constants $\eta_0$  and
$C$ depending only on $d$ and $p$ such that for every $\eta \leq
\eta_0$, we have 
\begin{equation}\label{packresulteq}
N(\eta; p, \Gamma)  \geq  \exp \left(C \left(
    \frac{\Gamma}{\eta}\right)^{(d-1)/2}\right). 
\end{equation}
\end{theorem}
\begin{IEEEproof}
Observe that by scaling, it is enough to prove for the case $\Gamma =
1$. We loosely follow Bronshtein~\cite[Proof of Theorem
4]{Bronshtein76}. Fix $\epsilon \in (0,
1)$. For each point $x \in S^{d-1}$, let $S_x$ denote the supporting
hyperplane to the unit ball $B$ at $x$ and let $H_x$ be the hyperplane
intersecting the sphere that is parallel to $S_x$ and at a distance of
$\epsilon$ from $S_x$. Let $H_x^+$ and $H_x^-$ denote the two
halfspaces bounded by $H_x$ where we assume that $H_x^+$ contains the
origin. Let $T_x := S^{d-1} \cap H_x^-$ and $A_x := B \cap H_x$, where
$B$ stands for the unit ball. It can be checked that the
(Euclidean) distance between $x$ and every point in $T_x$ (and $A_x$) is
less than or equal to $\sqrt{2}\sqrt{\epsilon}$. It follows that if
the distance between two points $x$ and $y$ in $S^{d-1}$ is strictly
larger than $2\sqrt{2}\sqrt{\epsilon}$, then the sets $T_x$ and $T_y$ are
disjoint. 

By standard results (see for example~\cite[Proof of Theorem
4]{Bronshtein76} where it is referred to as Mikhlin's result), there
exist positive constants $C_1$, depending 
only on $d$, and $\epsilon_0$ such that for every $\epsilon \leq
\epsilon_0$, there exist $N \geq C_1(\sqrt{\epsilon})^{1-d}$ points
$x_1, \dots, x_N$ in $S^{d-1}$ such that the Euclidean distance
between $x_i$ and $x_j$ is strictly larger than
$2\sqrt{2}\sqrt{\epsilon}$ whenever $i \neq j$. From now on, we assume that 
$\epsilon \leq \epsilon_0$. We then consider a mapping $\Phi : \left\{
0, 1\right\}^N \rightarrow \K^d(1)$, which is defined, for $\tau = (\tau_1,
\dots, \tau_N) \in \left\{0, 1 \right\}^N$, by  
\begin{equation*}
\Phi(\tau) := B \cap D_1(\tau_1) \cap D_2(\tau_2) \cap \dots \cap D_N(\tau_N),
\end{equation*}
where for $i = 1, \dots, N$,
\begin{equation*}
D_i(0)  :=  H_{x_i}^+  \text{ and }  D_i(1) := B. 
\end{equation*}
It must be clear that the Hausdorff distance between $\Phi(\tau)$ and
$\Phi(\tau ')$ is not less than $\epsilon$ (in fact, it is exactly equal
to $\epsilon$) if $\tau \neq \tau '$. Thus, $\left\{\Phi(\tau):
  \tau \in \left\{0, 1\right\}^N \right\}$ is an $\epsilon$-packing
set for $\K^d(1)$ under the Hausdorff metric. However, it is \textit{not} an
$\epsilon$-packing set under the $L^p$ metric. Indeed, the $L^p$
distance between $\Phi(\tau)$ and $\Phi(\tau ')$ is not necessarily larger than
$\epsilon$ for all pairs $(\tau, \tau '), \tau \neq \tau '$. The
$L^p$ distance between $\Phi(\tau)$ and $\Phi(\tau ')$ depends on
the Hamming distance $\Upsilon(\tau, \tau ') = \sum_i \left\{ \tau_i
  \neq \tau_i'\right\}$ between $\tau$ and $\tau '$. We make the claim that  
\begin{equation}\label{claim}
\delta_p\left(\Phi(\tau), \Phi(\tau ') \right)  \geq  C_2 \epsilon
\left( \sqrt{\epsilon}\right)^{(d-1)/p} \left( \Upsilon(\tau, \tau ')
\right)^{1/p}, 
\end{equation}
where $C_2$ depends only on $d$ and $p$.
The claim will be proved later. Assuming it is true, we recall the
Varshamov-Gilbert lemma from the previous section to assert the
existence of a subset $W$ of $\left\{0, 1\right\}^N$ with $|W| \geq
\exp(N/8)$ such that $\Upsilon(\tau, \tau') \geq N/4$ for all $\tau,
\tau' \in W$ with $\tau \neq \tau'$. Because 
$N \geq C_1(\sqrt{\epsilon})^{1-d}$, we get from~\eqref{claim} that
for all $\tau, \tau' \in W$ with $\tau \neq \tau'$, we have 
\begin{equation*}
\delta_p\left(\Phi(\tau), \Phi(\tau') \right)  \geq  C_3 \epsilon  \text{
  where }  C_3 := C_2 \left(
  \frac{C_1}{4}\right)^{1/p}. 
\end{equation*}
Taking $\eta := C_3 \epsilon$, we have obtained, for each $\eta \leq
\eta_0 := C_3 \epsilon_0$, an $\eta$-packing subset of $\K^d(1)$ with
size $M$, where 
\begin{equation*}
\log M  \geq  N/8  \geq  \frac{C_1}{8} \left(
  \frac{1}{\sqrt{\epsilon}}\right)^{d-1}  =  C_4 \left(
  \frac{1}{\sqrt{\eta}}\right)^{d-1}. 
\end{equation*}
The constant $C_4$ only depends on $d$ and $p$ thereby
proving~\eqref{packresulteq}.  

It remains to prove the claim~\eqref{claim}. Fix a point $x \in
S^{d-1}$ and $\epsilon \in (0,1)$. We first observe that it is enough
to prove that  
\begin{equation}\label{tworeg}
 \delta_p(A_x, T_x)^p  \geq  C_5 \epsilon^p
\left( \sqrt{\epsilon}\right)^{d-1}, 
\end{equation} 
for a constant $C_5$ depending on just $d$ and $p$, where $A_x$ and
$T_x$ are as defined in the beginning of the proof. This is because of
the fact that for every $\tau, \tau' \in W$ with $\tau \neq \tau'$, we
can write  
\begin{equation}\label{addit}
\delta_p\left(\Phi(\tau),\Phi(\tau')\right)^p  =
\sum_{i \in I} \delta_p(A_{x_i}, T_{x_i})^p, 
\end{equation}
where $I := \left\{1\leq i \leq N : \tau_i \neq \tau_i'\right\}$. The
equality~\eqref{addit} is a consequence of the fact that the points
$x_1, \dots, x_N$ are chosen so that $T_{x_1}, \dots, T_{x_N}$ are
disjoint. 

We shall now prove the inequality~\eqref{tworeg} which will complete
the proof. Let $u_0$ denote the point in $A_x$ that is closest to the
origin. Also let $u_1$ be a point in $A_x \cap
S^{d-1}$. Let $\alpha$ denote the angle between $u_0$ and
$u_1$. Clearly, $\alpha$ does not depend on the choice of $u_1$ and
$\cos \alpha = 1 - \epsilon$. Now let $u$ be a fixed unit vector and
let $\theta$ be the angle between the vectors $u$ and $u_0$. By
elementary geometry, we deduce that 
\begin{equation*}
h_{T_x}(u) - h_{A_x}(u)  =  
\begin{cases} 1 -  \cos \left(\alpha - \theta \right) & \text{if $0
    \leq \theta \leq \alpha$,} 
\\
0 &\text{otherwise.}
\end{cases}
\end{equation*}
Because the difference of support functions only depends on the angle
$\theta$, we can write, for a constant $C_6$ depending only on $d$,
that  
\begin{equation*}
\delta_p(A_x, T_x)^p  =  C_6 \int_0^{\alpha}
\left( 1 - \cos (\alpha - \theta)\right)^p \sin^{d-2} \theta d\theta.
\end{equation*}
Now suppose $\beta$ is such that $\cos (\alpha - \beta) = 1 -
\epsilon/2$. Then from above, we get that 
\begin{align*}
\delta_p(A_x, T_x)^p  &\geq  C_6
\int_0^{\beta} \left( 1 - \cos (\alpha - \theta)\right)^p \sin^{d-2}\theta d\theta \\ 
&\geq  C_6 \left(\frac{\epsilon}{2}\right)^p \int_0^{\beta}
\sin^{d-2}\theta  d\theta \\ 
&\geq  C_6 \left(\frac{\epsilon}{2}\right)^p \int_0^{\beta}
\sin^{d-2} \theta \cos \theta d\theta \\ 
&=  \frac{C_6}{d-1} \left(\frac{\epsilon}{2}\right)^p  \sin^{d-1}
  \beta.
\end{align*}
We shall show that $\sin \beta \geq
\left(\sqrt{\epsilon}\right)/(2\sqrt{2})$ which will
prove~\eqref{tworeg}. Recall that $\cos \alpha = 1 - \epsilon$. Thus 
\begin{align*}
1 - \frac{\epsilon}{2}  &= \cos (\alpha - \beta) \\ 
&\leq  \cos \alpha + \sin \alpha \sin \beta \\
&=  1 - \epsilon + \sqrt{1 - (1-\epsilon)^2} \sin \beta \\
&\leq  1 - \epsilon + \sqrt{2} \sqrt{\epsilon} \sin \beta,
\end{align*}
which when rearranged gives $\sin \beta \geq
\left(\sqrt{\epsilon}\right)/(2\sqrt{2})$. The proof is complete. 
\end{IEEEproof}

\section{Conclusion}
By a simple application of convexity, we proved an inequality relating
the minimax risk in multiple hypothesis testing problems to
$f$-divergences of the probability measures involved. This inequality
is an extension of Fano's inequality. As another corollary, we obtained a sharp
inequality between total variation distance and $f$-divergences. We
also indicated how to control the quantity $J_f$ which appears in our
lower bounds. This leads to important global lower bounds for the minimax
risk. Two applications of our bounds are presented. In the first
application, we used the bound~\eqref{myb.chi} to prove a new lower
bound (which turns to be rate-optimal) for the minimax risk of
estimating a convex body from noisy measurements of the support
function in $n$ directions. In the second application, we employed
inequalities~\eqref{normfano} and~\eqref{jsyb} to give a
different proof of a recent lower bound for covariance matrix
estimation due to Cai, Zhang and Zhou~\cite{CaiZhangZhou2009}.


%



\section*{Acknowledgment}
The author is indebted to David Pollard for his insight and also for
numerous stimulating discussions which led to many of the ideas in
this paper; to Andrew Barron for his constant encouragement and
willingness to discuss his own work on minimax bounds. Thanks are also 
due to Aditya Mahajan for pointing out to the author
that inequality~\eqref{maha.eq1} has the extension~\eqref{maha.eq} for
the case of non-uniform priors $w$; to an anonymous referee for
helpful comments, for pointing out an error and for
Example~\ref{ddimnorm}; to Richard Gardner for comments that
greatly improved the quality of the paper and to Alexander Gushchin
for informing the author about his paper~\cite{Gushchin} and for
sending him a scanned copy of it.

\ifCLASSOPTIONcaptionsoff
  \newpage
\fi



%
\bibliographystyle{IEEEtran}
\def\noopsort#1{}

%

\begin{IEEEbiographynophoto}{Adityanand Guntuboyina}
(S'10) received the B.Stat. and M.Stat. degrees
from the Indian Statistical Institute, Kolkata, in 2004 and 2006
respectively and the M.A. degree in statistics from Yale University in
2008. 

He is currently a PhD candidate in the Department of Statistics at
Yale University. His dissertation advisor is Professor David
Pollard. His main research interests include various areas of
nonparametric statistics like high-dimensional nonparametric
regression and bayesian nonparametrics. He is also interested in
the concentration of measure phenomena, random matrices and respondent
driven sampling. 
\end{IEEEbiographynophoto}




\vfill


\end{document}